\input amssym.def
\input amssym.tex

\headline={\ifnum \pageno=1 {\hfill} 
\else{\hss \tenrm -- \folio\ -- \hss}\fi}
\footline={\hfil}

\def\dater{\vglue-10mm\rightline{(\the\day/\the\month/\the\year)}}

\hsize 146mm
\vsize 224mm
\hoffset=6mm
\voffset=8mm
\baselineskip=5mm
\overfullrule =0pt

\font\ninerm=cmr9
 at 10,5pt
\font\pptitre=cmbx10 at 11pt

\font\Gtitre=cmbx10 at 15pt


\def\og{\leavevmode\raise.30ex
\hbox{$\scriptscriptstyle\langle\!\langle\>$}}    
\def\bigog{\leavevmode\raise.30ex
\hbox{$\langle\!\langle\>$}}  
\def\fg{\leavevmode\raise.30ex
\hbox{$\scriptscriptstyle\>\rangle\!\rangle$}}    
\def\bigfg{\leavevmode\raise.30ex
\hbox{$\>\rangle\!\rangle$}}    

\catcode`\@=11

\font\author=cmcsc10
\font\pauthor=cmcsc10 at 8pt
\font\tenmsx=msam10
\font\sevenmsx=msam10 scaled 700
\font\fivemsx=msam10 scaled 500
\font\tenmsy=msbm10
\font\sevenmsy=msbm10 scaled 700
\font\fivemsy=msbm10 scaled 500
\newfam\msxfam
\newfam\msyfam
\textfont\msxfam=\tenmsx  \scriptfont\msxfam=\sevenmsx
\scriptscriptfont\msxfam=\fivemsx
\textfont\msyfam=\tenmsy  \scriptfont\msyfam=\sevenmsy
\scriptscriptfont\msyfam=\fivemsy

\def\hexnumber@#1{\ifnum#1<10 \number#1\else
\ifnum#1=10 A\else\ifnum#1=11 B\else\ifnum#1=12 C\else
\ifnum#1=13 D\else\ifnum#1=14 E\else\ifnum#1=15 F\fi\fi\fi\fi\fi\fi\fi}

\def\msx@{\hexnumber@\msxfam}
\def\msy@{\hexnumber@\msyfam}
\mathchardef\nmid="3\msy@2D
\mathchardef\varnothing="0\msy@3F
\mathchardef\nexists="0\msy@40
\mathchardef\smallsetminus="2\msy@72
\def\Bbb{\ifmmode\let\next\Bbb@\else
\def\next{\errmessage{Use \string\Bbb\space only in math mode}}\fi\next}
\def\Bbb@#1{{\Bbb@@{#1}}}
\def\Bbb@@#1{\fam\msyfam#1}

\font\tentbl=cmr10 scaled 900
\font\seventbl=cmr7 scaled 900
\font\fivetbl=cmr5 scaled 900

\newfam\tblfam

\textfont\tblfam=\tentbl
\scriptfont\tblfam=\seventbl
\scriptscriptfont\tblfam=\fivetbl


\def \C {{\Bbb C}}

\def \E {{\Bbb E}}

\def \N {{\Bbb N}}

\def \R {{\Bbb R}}

\font\ccm=cmmi10 at 9pt
\def\ccmL{\hbox{\ccm L}}

 at 9,5pt

\font\cmm=cmmi10 at 15pt
\def\cmmL{\hbox{\cmm L}}

\def\cmmL{\hbox{\cmm L}}

\font\f=cmr10 at 7pt

\def\f1{\hbox{\f 1}}
\def\f2{\hbox{\f 2}}
\def\f3{\hbox{\f 3}}
\def\f4{\hbox{\f 4}}
\def\f5{\hbox{\f 5}}
\def\f6{\hbox{\f 6}}
\def\f7{\hbox{\f 7}}
\def\f8{\hbox{\f 8}}
\def\f9{\hbox{\f 9}}

\def \d {\,{\rm d}}
\def\re{{\Re e\,}}

\def \dm {{\hbox {${1\over 2}$}}}
\def \sset {{\smallsetminus }}

\def\le{\leqslant}
\def\ge{\geqslant}

\topskip=12pt
\font\sept=cmti9

\def\rightheadline{\ifnum\pageno=\chstart{\hfill}
         \else\centerline{\sept On a conjecture of Montgomery-Vaughan
on extreme values of automorphic $\ccmL$-functions at 1}\hfill
\hskip -3,5mm \tenrm\folio\fi}
\def\leftheadline{\ifnum\pageno=\chstart{\hfill}
         \else\tenrm\folio \hskip -3,5mm \hfill
\centerline{\pauthor J.-Y. Liu, E. Royer \& J. Wu}\fi}
\headline={\ifnum\pageno=\chstart{\hfill}
\else{\ifodd\pageno\rightheadline\else\leftheadline\fi}\fi}
\footline={\hfill}

\pageno=1
\newcount\chstart
\chstart=\pageno


\vglue 5mm

\centerline{\Gtitre On a conjecture of Montgomery-Vaughan}

\bigskip

\centerline{\Gtitre 
on extreme values of automorphic $\cmmL$-functions at 1}

\vskip 8mm

\centerline{\author J.-Y. Liu, E. Royer \& J. Wu}

\vskip 10mm

{\leftskip=1cm
\rightskip=1cm
{\ninerm
{\bf Abstract}.
In this paper, 
we prove a weaker form of a conjecture of Montgomery-Vaughan
on extreme values of automorphic $\ccmL$-functions at 1.
\par}}

\vskip 10mm

{\hsize 122mm
\baselineskip=4mm
\overfullrule =0pt
\leftskip=18mm

\hskip 40mm{\bf Contents}

\vskip 5mm

\ninerm
\S\ 1. Introduction
\dotfill \hskip 6mm 1

\vskip 0,3mm

\S\ 2. Expression of $E(s, y)$ and existence of saddle-point
\dotfill \hskip 6mm 7

\vskip 0,3mm

\S\ 3. Preliminary lemmas 
\dotfill \hskip 6mm 9

\vskip 0,3mm

\S\ 4. Estimates of $\phi_n(\sigma, y)$ 
\dotfill \hskip 6mm 11

\vskip 0,3mm

\S\ 5. Estimates of $|E(\sigma+i\tau, y)|$ 
\dotfill \hskip 6mm 18

\vskip 0,3mm

\S\ 6. Proof of Theorem 3
\dotfill \hskip 6mm 21

\vskip 0,3mm

\S\ 7. Proof of Theorem 4
\dotfill \hskip 6mm 24

\vskip 0,3mm

\S\ 8. Proof of Corollary 5
\dotfill \hskip 6mm 24

\vskip 0,3mm

\S\ 9. Proof of Theorem 2
\dotfill \hskip 6mm 26

\vskip 0,3mm

\S\ 10. Proof of Theorem 1
\dotfill \hskip 6mm 30

References
\dotfill \hskip 4,34mm 30
\par}

\vskip 10mm

\centerline{\pptitre \S\ 1. Introduction}

\bigskip

The automorphic $L$-functions constitute a powerful tool 
for studying arithmetic,  algebraic or geometric objects.
For squarefree integer $N$ 
and even integer $k$,
denote by ${\rm H}_k^*(N)$ the set of all newforms of level $N$
and of weight $k$. 
It is known that 
$$|{\rm H}_k^*(N)|
= {k-1\over 12}\varphi(N)
+ O\big((kN)^{2/3}\big),
\leqno(1.1)$$
where $\varphi(N)$ is the Euler function
and the implied constant is absolute.
Let $m\ge 1$ be an integer and
let $L(s, {\rm sym}^mf)$ be the $m$th symmetric power $L$-function 
of $f\in {\rm H}_k^*(N)$
normalised so that the critical strip is given by $0<\re s<1$. 
The values of these functions at the edge of the critical strip
contain information of great interest.
For example, 
Serre [18] showed that
the Sato-Tate conjecture is equivalent to 
$L(1+i\tau, {\rm sym}^mf)\not=0$
for all $m\in \N$ and $\tau\in \R$.
The distribution of the values $L(1, {\rm sym}^mf)$
has received attention of many authors, including
Goldfeld, Hoffstein \& Lieman [2],
Hoffstein \& Lockhart [7],
Luo [12], 
Royer [14, 15], 
Royer \& Wu [16, 17], 
Cogdell \& Michel [1], 
Habsieger \& Royer [5]
and Lau \& Wu [10, 11].
In particular, Lau \& Wu ([10], [11])
proved the following results:

(i)
For every fixed integer $m\ge 1$,
there are four positive constants $A_m^\pm$ and $B_m^\pm$ such that
for any newform $f\in {\rm H}_k^*(1)$,
under the Great Riemann Hypothesis (GRH) for $L(s, {\rm sym}^mf)$, 
we have, for $k\to\infty$, 
$$\{1+o(1)\} (2B_m^-\log_2k)^{-A_m^-}
\le L(1, {\rm sym}^mf)
\le \{1 + o(1)\} (2B_m^+\log_2k)^{A_m^+}.
\leqno(1.2)$$ 
Here (and in the sequel) 
$\log_j$ denotes the $j$-fold iterated logarithm.
For most values of $m$,
the constants $A_m^\pm$ and $B_m^\pm$ 
can be explicitly evaluated, for example,
$$\cases{
A_m^+ = m+1,
\qquad
B_m^+ = e^\gamma
\hskip 18,5mm
(m\in \N),
\cr\noalign{\smallskip}
A_m^- = m+1,
\qquad
B_m^- = e^\gamma \zeta(2)^{-1}
\hskip 8,5mm
({\rm odd}\; m),
\cr\noalign{\smallskip}
A_2^- \hskip 0,3mm = 1,
\hskip 14,5mm
B_2^- = e^\gamma \zeta(2)^{-2},
\cr\noalign{\smallskip}
A_4^- = {5\over 4},
\hskip 14,5mm
B_4^- = e^\gamma B_4'^{-},
\cr}$$
where $\zeta(s)$ is the Riemann zeta-function,
$\gamma$ denotes the Euler constant and  
$B_4'^{-}$ is a positive constant given 
by a rather complicated 
Euler product ([10], Theorem 3).

(ii)
In the opposite direction, it was shown unconditionally
that for $m\in \{1, 2, 3, 4\}$ 
there are newforms $f_m^{\pm}\in {\rm H}_k^*(1)$
such that for $k\to\infty$ ([10], Theorem 2),
$$\cases{
L(1, {\rm sym}^mf_m^+)
\ge \{1 + o(1)\} (B_m^+ \log_2k)^{A_m^+},
\cr\noalign{\medskip}
L(1, {\rm sym}^mf_m^-)
\le \{1+o(1)\} (B_m^-\log_2k)^{-A_m^-}.
\cr}
\leqno(1.3)$$

(iii)
In the aim of removing GRH and 
closing up the gap coming from the factor 2 in (1.2)
(comparing it with (1.3)),
an almost all result was established.
Let $\varepsilon>0$ be an arbitrarily small positive number,
$m\in \{1, 2, 3, 4\}$ and $2\mid k$.
Then there is a subset ${\rm E}_k^*$ of ${\rm H}_k^*(1)$ 
such that
$|{\rm E}_k^*|\ll {\rm H}_k^*(1) e^{-(\log k)^{1/2-\varepsilon}}$
and for each $f\in {\rm H}_k^*(1)\sset {\rm E}_k^*$, 
we have, for $k\to\infty$,
$$\{1+O(\varepsilon_k)\} (B_m^-\log_2k)^{-A_m^-}
\le L(1, {\rm sym}^mf)
\le \{1+O(\varepsilon_k)\} (B_m^+\log_2k)^{A_m^+},
\leqno(1.4)$$
where $\varepsilon_k:=(\log k)^{-\varepsilon}$
and the implied constants depend on $\varepsilon$ only
([11], Corollary 2).

By comparing (1.3) with (1.4),
the extreme values of $L(1, {\rm sym}^mf)$
seem to be  given by (1.3).
Clearly it is interesting to investigate further 
the size of exceptional set ${\rm E}_k^*$.
In the case of quadratic characters $L$-functions,
Montgomery \& Vaughan [13] proposed, 
based on a probabilistic model,
three conjectures on the size of exceptional set.
The first one has been proved recently 
by Granville \& Soundararajan [4].
As Cogdell \& Michel indicated in [1],
it would be interesting to try to get,
as close as possible, the analogues of the conjectures of 
Montgomery-Vaughan for automorphic $L$-functions.
The analogue of Montgomery-Vaughan's first conjecture
for the automorphic symmetric power $L$-functions 
can be stated as follows.

\proclaim Conjecture.
Let $m\ge 1$ be a fixed integer and
$$\eqalign{
F_k(t, {\rm sym}^m)
& := {1\over |{\rm H}_k^*(1)|}
\sum_{f\in {\rm H}_k^*(1),\; L(1, {\rm sym}^mf)\ge (B_m^+t)^{A_m^+}} 1,
\cr\noalign{\vskip 2mm}
G_k(t, {\rm sym}^m)
& := {1\over |{\rm H}_k^*(1)|}
\sum_{f\in {\rm H}_k^*(1),\; 
L(1, {\rm sym}^mf)\le (B_m^-t)^{-A_m^-}} 1.
\cr}$$
Then there are positive constants $c_i=c_i(m)\;(i=1,2)$ such that
for $k\to\infty$,
$$\cases{
e^{-c_1(\log k)/\log_2k}
\ll F_k(\log_2k, {\rm sym}^m)
\ll e^{-c_2(\log k)/\log_2k},
\cr\noalign{\vskip 3,5mm}
e^{-c_1(\log k)/\log_2k}
\ll G_k(\log_2k, {\rm sym}^m)
\ll e^{-c_2(\log k)/\log_2k}.
\cr}
\leqno(1.5)$$

The aim of this paper is 
to prove a weaker form of this conjecture for $m=1$.
In this case, we write, for simplification of notation, 
$$L(s, f)=L(s, {\rm sym}^1f),
\qquad
F_k(t)=F_k(t, {\rm sym}^1),
\qquad
G_k(t)=G_k(t, {\rm sym}^1).$$
In view of the trace formula of Petersson ([8], Theorem 3.6),
it is more convenient to consider the weighted arithmetic distribution
function.
As usual, denote by 
$$\omega_f 
:= {\Gamma(k-1)\over (4\pi)^{k-1}\|f\|}$$
the harmonic weight in modular forms theory
and define the weighted arithmetic distribution functions
$$\eqalign{\widetilde{F}_k(t)
& := \Big(\sum_{f\in {\rm H}_k^*(1)} \omega_f\Big)^{-1}
\sum_{f\in {\rm H}_k^*(1), \;
L(1, f)\ge (e^\gamma t)^2} \omega_f,
\cr\noalign{\vskip 3mm}
\widetilde{G}_k(t)
& := \Big(\sum_{f\in {\rm H}_k^*(1)} \omega_f\Big)^{-1}
\sum_{f\in {\rm H}_k^*(1), \;
L(1, f)\le (6\pi^{-2}e^\gamma t)^{-2}} \omega_f.
\cr}$$
By using (1.1),  the classical estimate
$$\mathop{\sum_{f\in {\rm H}_k^*(1)}} \omega_f
= 1+O\big(k^{-5/6}\big)
\leqno(1.6)$$
and the bound of Goldfeld, Hoffstein \& Lieman [2]: 
$$1/(k\log k)\ll \omega_f \ll (\log k)/k,
\leqno(1.7)$$
we easily see that 
$$\cases{
\widetilde{F}_k(t)/\log k 
\ll F_k(t)
\ll \widetilde{F}_k(t)\log k,
\cr\noalign{\vskip 2mm}
\widetilde{G}_k(t)/\log k 
\ll G_k(t)
\ll \widetilde{G}_k(t)\log k.
\cr}
\leqno(1.8)$$
This shows that in order to prove (1.5)
it is sufficient to establish corresponding estimates 
of the same quality for $\widetilde{F}_k(t)$
and $\widetilde{G}_k(t)$.

Our main result is the following one.

\proclaim Theorem 1.
For any $A\ge 1$ there are two positive constants
$c=c(A)$ and $C=C(A)$ such that the estimate
$$\widetilde{F}_k(t)
= \{1+\Delta_k(t)\}
\exp\bigg\{-{e^{t-\gamma_0}\over t}
\bigg(1+O\bigg({1\over t}\bigg)\bigg)\bigg\}
\leqno(1.9)$$
holds uniformly for $k\ge 16, 2\mid k$ and $t\le T(k)$,
where $\gamma_0$ is given by (1.24) below,
$|\theta|\le 1$ and
$$\cases{
\Delta_k(t):=\theta e^{t-T(k)-C}(t/T(k))^{1/2}
+ O_A\big(e^{-ce^{t/5}}+(\log k)^{-A}\big),
\cr\noalign{\vskip 2mm}
T(k):=\log_2k-\textstyle {5\over 2}\log_3k-\log_4k-3C.
\cr}
\leqno(1.10)$$
In particular 
there are two positive constants $c_1$ and $c_2$ such that
$$e^{-c_1(\log k)/\{(\log_2k)^{7/2}\log_3k\}}
\ll F_k(T(k))
\ll e^{-c_2(\log k)/\{(\log_2k)^{7/2}\log_3k\}}.
\leqno(1.11)$$
The similar estimates for $\widetilde{G}_k(t)$
and $G_k(T(k))$ hold also.

\goodbreak
\smallskip

{\bf Remark 1}.
The estimates (1.11) of Theorem 1 can be considered as
a weaker form of Montgomery-Vaughan's conjecture (1.5) for $m=1$,
since $T(k)\sim \log_2k$ as $k\to\infty$.
Moreover, if we could take $T(k)=\log_2k$ in (1.11) then (1.9) would lead
to the Montgomery-Vaughan's conjecture (1.5). Hence we fail from a
shift 
$$
{5\over 2}\log_3k+\log_4k+3C.
$$
It seems however to be rather difficult to resolve completely this conjecture.
One of the main difficulties is that there are no analogues
of the quadratic reciprocity law and 
Graham-Ringrose's estimates for short characters sums 
of friable moduli [3],
which have been exploited by Granville \& Soundararajan [4].

\medskip

In order to prove Theorem 1,
we need to introduce a probabilistic model as in [1].
Consider a probability space $(\Omega, \mu)$, with measure $\mu$. 
Let $\hbox{SU(2)}^\natural$ be 
the set of conjugacy classes of $\hbox{SU(2)}$. 
The group $\hbox{SU(2)}$ is endowed 
with its Haar measure $\mu_{\rm H}$ and
$$
\hbox{SU(2)}^\natural
= \bigg\{
\pmatrix{e^{i\theta} & 0 \cr 0 & e^{-i\theta}\cr} : \, \theta\in[0,\pi]
\bigg\}\biggr/\sim\biggr.
$$
is endowed with the Sato-Tate measure 
$\d \mu_{\rm st}(\theta) : = (2/\pi)\sin^2\theta \d \theta$,
{\it i.e.,\null}
the direct image of $\mu_{\rm H}$ 
by the canonical projection $\hbox{SU(2)}\to\hbox{SU(2)}^\natural$.
On the space $(\Omega, \mu)$, 
define a sequence indexed by the prime numbers, 
$g^\natural(\omega)=\{g^\natural_p(\omega)\}_p$ of random matrices
taking values in $\hbox{SU(2)}^\natural$, 
given by
$$
g_p^\natural(\omega)
:=\pmatrix{
e^{i\vartheta_p(\omega)} & 0
\cr\noalign{\vskip 1mm}
0                        & e^{-i\vartheta_p(\omega)}
\cr}^\natural.
$$
We assume that each function $g^\natural_p(\omega)$ is distributed 
according to the Sato-Tate measure. 
This means that, for each integrable function 
$\phi : \hbox{SU(2)}^\natural \to \R$, 
the expected value of $\phi\circ g^\natural_p$ is
$$
\E(\phi\circ g^\natural_p)
:=\int_\Omega\phi\circ g^\natural_p(\omega)\d\mu(\omega)
=\int_0^\pi\phi\biggl(
\pmatrix{e^{i\theta} & 0 \cr 0 & e^{-i\theta}\cr}
\biggr)\cdot (2/\pi)\sin^2\theta \d \theta.
$$
Moreover, we assume that the sequence $g^\natural(\omega)$ is made of
independent random variables. 
This means that, for any sequence of integrable functions
$\{G_p : \hbox{SU(2)}^\natural \to \R\}_p$, we have
$$\leqalignno{
\E\Big(\prod_p G_p\circ g^\natural_p\Big)
& := \int_\Omega\prod_p G_p\circ g^\natural_p(\omega)\d\mu(\omega)
& (1.12)
\cr
& \, = \prod_p\int_\Omega G_p\circ g^\natural_p(\omega)\d\mu(\omega)
&
\cr
& \, = \prod_p\int_0^\pi
G_p\biggl(
\pmatrix{e^{i\theta} & 0 \cr 0 & e^{-i\theta}\cr}
\biggr)\cdot (2/\pi)\sin^2\theta \d \theta.
\cr}$$
Let $I$ be the identity matrix. Then for $\re s>\dm$, 
the random Euler product
$$L(s, g^\natural(\omega))
:= \prod_p
\det\big(I-p^{-s}g_p^\natural(\omega)\big)^{-1}
=:\prod_p L_p(s, g^\natural(\omega))$$
turns out to be absolutely convergent a.s.

Now we define our probabilistic distribution functions
$$\cases{
\Phi(t)
:= {\rm Prob}
\big(\big\{L(1, g^\natural(\cdot))
\ge (e^\gamma t)^2\big\}\big),
\cr\noalign{\vskip 2mm}
\Psi(t)
:= {\rm Prob}\big(
\big\{L(1, g^\natural(\cdot))\le (6\pi^{-2}e^\gamma t)^{-2}\big\}
\big).
\cr}$$

We shall prove Theorem 1 in two steps.
The first one is to compare 
$\widetilde{F}_k(t)$ with $\Phi(t)$ 
(resp. $\widetilde{G}_k(t)$ with $\Psi(t)$). 

\proclaim Theorem 2.
For any $A\ge 1$ there are two positive constants
$c=c(A)$ and $C=C(A)$ such that
the asymptotic formulas
$$\widetilde{F}_k(t)
=\Phi(t)\{1+\Delta_k(t)\}
\qquad
\hbox{and} 
\qquad
\widetilde{G}_k(t)
=\Psi(t)\{1+\Delta_k(t)\}
\leqno(1.13)$$
hold uniformly for $k\ge 16, 2\mid k$ and $t\le T(k)$,
where $\Delta_k(t)$ and $T(k)$ are defined  by (1.10).

\medskip

The second step of the proof of Theorem 1 is the evaluation of $\Phi(t)$ (resp. $\Psi(t)$).
For this, 
we consider a truncated random Euler product
$$L(s, g^\natural(\omega); y)
:= \prod_{p\le y} L_p(s, g^\natural(\omega))$$
and the corresponding distribution functions
$$\cases{
\Phi(t, y)
:= {\rm Prob}
\big(\big\{L(1, g^\natural(\omega); y)
\ge (e^\gamma t)^2\big\}\big),
\cr\noalign{\vskip 2mm}
\Psi(t, y)
:= {\rm Prob}\big(
\big\{L(1, g^\natural(\omega); y)\le (6\pi^{-2}e^\gamma t)^{-2}\big\}
\big).
\cr}$$
We have
$$\Phi(t) = \Phi(t, \infty)
\qquad{\rm and}\qquad
\Psi(t) = \Psi(t, \infty).
\leqno(1.14)$$

We shall use the saddle-point method 
(introducted by Hildebrand \& Tenenbaum [6])
to evaluate $\Phi(t, y)$ and  $\Psi(t, y)$.
For this, we need to introduce some notation.
For $s\in \C$ and $y\ge 2$, define
$$E(s, y) := \E\big(L(1, g^\natural(\omega); y)^s\big)
\qquad{\rm and}\qquad
E(s) := E(s, \infty),
\leqno(1.15)$$
where $\E(\cdot)$ denotes the expected value.
We define also
$$\phi(s, y) := \log E(s, y),
\qquad
\phi_n(s, y) 
:= \displaystyle {\partial^n\phi\over \partial s^n}(s, y)
\quad(n\ge 0).
\leqno(1.16)$$

According to Lemmas 2.3 and 8.1 below,
there is an absolute constant $c\ge 2$ such that
for $t\ge 4\log c$ and $y\ge ce^t$, 
the equation
$$\phi_1(\kappa, y) = 2(\log t + \gamma)
\leqno(1.17)$$
has a unique positive solution $\kappa=\kappa(t, y)$
and for each integer $J\ge 1$, 
there are computable constants 
$\gamma_0, \gamma_1, \dots, \gamma_J$ such that
the asymptotic formula
$$\kappa(t, y)
= e^{t-\gamma_0}
\bigg\{1 + \sum_{j=1}^J {\gamma_j\over t^j} 
+ O_J\bigg({1\over t^{J+1}}
+{e^tt\over y\log y}\bigg)
\bigg\}
\leqno(1.18)$$
holds uniformly for $t\ge 1$ and $y\ge 2e^t$,
the constant $\gamma_0$ beign given by (1.24) below.

Finally write $\sigma_n := \phi_n(\kappa, y)$.

\smallskip

\proclaim Theorem 3.
We have
$$
\Phi(t, y)
= {E(\kappa, y)\over \kappa\sqrt{2\pi\sigma_2}(e^\gamma t)^{2\kappa}}
\bigg\{1 + O\bigg({t\over e^t}\bigg)\bigg\}
$$
uniformly for $t\ge 1$ and $y\ge 2e^t$.

\proclaim Theorem 4.
For each integer $J\ge 1$, we have
$$\Phi(t, y)
=  \exp\bigg\{
- \kappa
\bigg[
\sum_{j=1}^J {a_j\over (\log\kappa)^j}
+ O_J\big(R_J(\kappa, y)\big)
\bigg]\bigg\}
\leqno(1.19)$$
uniformly for $t\ge 1$ and $y\ge 2e^t$, 
where the error term $R_J(\kappa, y)$ is given by 
$$R_J(\kappa, y)
:= {1\over (\log\kappa)^{J+1}}
+ {\kappa\over y\log y}
\leqno(1.20)$$
and
$$a_j :=  
\int_0^\infty \bigg({h(u)\over u}\bigg)' (\log u)^{j-1} \d u
\leqno(1.21)$$
with
$$h(u) := \cases{
\displaystyle
\log\bigg(
{2\over \pi} \int_0^{\pi} 
e^{2u\cos\theta}
\sin^2\theta \d\theta
\bigg)
& if $\;\;0\le u<1$,
\cr\noalign{\medskip}
\displaystyle
\log\bigg(
{2\over \pi} \int_0^{\pi} 
e^{2u\cos\theta}
\sin^2\theta \d\theta
\bigg)
- 2u
& if $\;\;u\ge 1$.
\cr}
\leqno(1.22)$$

\smallskip

As a corollary of Theorem 4,
we can obtain an asymptotic developpment
for $\log\Phi(t, y)$ in $t^{-1}$.
In particular we see that the probabilistic distribution function 
$\Phi(t)$ decays double exponentially as $t\to\infty$.

\proclaim Corollary 5.
For each integer $J\ge 1$, 
there are computable constants $a^*_1, \dots, a^*_J$ such that
the asymptotic formula
$$\Phi(t, y)
= \exp\bigg\{-e^{t-\gamma_0}
\bigg[\sum_{j=1}^J {a^*_j\over t^j}
+ O_J\big(R_J(e^t, y)\big)\bigg]\bigg\}
\leqno(1.23)$$
holds uniformly for $t\ge 1$ and $y\ge 2e^t$.
Further  we have
$$
\gamma_0 := {1\over 2}\int_0^\infty {h'(u)\over u} \d u, 
\qquad
a^*_1 : = 1,
\qquad
a^*_2 := \gamma_0 
- {\gamma_0^2\over 2}
- \int_0^\infty {h(u)\over u^2} (\log u) \d u.
\leqno(1.24)$$
In particular for each integer $J\ge 1$, we have
$$\Phi(t)
= \exp\bigg\{-e^{t-\gamma_0}
\bigg[\sum_{j=1}^J {a^*_j\over t^j}
+ O_J\bigg({1\over t^{J+1}}\bigg)\bigg]\bigg\}
\leqno(1.25)$$
uniformly for $t\ge 1$.

\goodbreak
\smallskip

{\bf Remark 2}.
(i)
The same results hold also for $\Psi(t, y)$.

(ii)
Taking $t=\log_2k$ and $J=1$ in (1.25) of Corollary 5, 
we see that the probabilistic distribution function $\Phi(t)$
(resp. $\Psi(t)$) verifies Montgomery-Vaughan's conjecture (1.5).
But (1.13) is too weak to derive this conjecture for 
$F_k(t)$ (resp. $G_k(t)$).
This means that we must take $T(k)=\log_2k$ in Theorem 2,
which seems be rather difficult.

(iii)
Our method can be generalized (with a little extra effort)
to prove that Theorems 1 and 2 hold for $L(1, {\rm sym}^mf)$ 
for $m\ge 1$ (unconditionally when $m=1, 2, 3, 4$ and 
under Cogdell-Michel's hypothesis ${\rm Sym}^m(f)$ and ${\rm LSZ}^m(1)$ 
[1] when $m\geqslant 5$)
and that Theorems 3, 4 and Corollary 5 are true for 
$L(1, {\rm sym}^mg^\natural(\omega); y)$ when $m\ge 1$. 
 
\medskip

\noindent{\bf Acknowledgements}.
We began working on this paper in November 2004 
during the visit of the first author 
to l'Institut Elie Cartan de Nancy,
and finished in January 2006 when the third author visited
School of Mathematics and System Sciences of Shandong University.
We are indebted to both institutions for invitations and support.
The second and third authors want to thank the CRM at Montr\'eal 
for its invitation. 
Finally we would express our sincere gratitude to Y.-K. Lau
of the University of Hong Kong for valuable discussion.

\vskip 8mm

\centerline{\pptitre \S\ 2. 
Expression of $E(s, y)$ and existence of saddle-point}

\bigskip

The aim of this section is to prove the existence of the saddle-point 
$\kappa(t, y)$, defined by equation (1.17).
The first step is to give an explicite expression of $E(s, y)$,
which is (1.24) of [1].
For the convenience of readers, 
we state it here as a lemma.

\proclaim Lemma 2.1.
For prime $p$, real $\theta$ and complex number $s$,
we define
$$D_p(\theta)
:= \prod_{0\le j\le 1} 
\big(1-e^{i(1-2j)\theta}p^{-1}\big)^{-1}
\qquad\hbox{and}\qquad
E_p(s) 
:= {2\over \pi} 
\int_0^\pi D_{p}(\theta)^{s} \sin^2\theta \d\theta.
\leqno(2.1)$$
Then for all $s\in \C$ and $y\ge 2$, we have
$$E(s, y) 
= \prod_{p\le y} E_p(s).
\leqno(2.2)$$

\noindent{\sl Proof}.
Taking
$$
G_p(M^\natural)=\cases{
\det\big(I-p^{-s'}M^\natural\big)^{-s} & if $\;p\le y$
\cr\noalign{\vskip 1mm}
1                                      & otherwise
\cr}$$
in (1.12), we get 
$$\eqalign{\E\big(L(s', g^\natural(\omega); y)^s\big)
& = \prod_{p\le y} 
\E\big(L_p(s', g_p^\natural(\omega))^s\big)
\cr
& = \prod_{p\le y} 
\int_{\Omega} 
\det\big(1-p^{-s'}g_p^\natural(\omega)\big)^{-s}
\d \mu(\omega)
\cr
& = \prod_{p\le y}
{2\over \pi} \int_0^{\pi}
\big(1 - 2p^{-s'}\cos\theta + p^{-2s'}\big)^{-s} 
\sin^2\theta \d\theta.
\cr}$$
Taking $s'=1$ and noticing (1.15) and (2.1),
we get the desired result.
\hfill
$\square$

\proclaim Lemma 2.2.
For all $p$ and $\sigma>0$, we have
$$E''_p(\sigma)E_p(\sigma)-E'_p(\sigma)^2>0.$$
In particular for all $\sigma>0$ and $y\ge 2$, 
we have $\phi_2(\sigma, y)>0$.

\noindent{\sl Proof}.
By using the definition (2.1) of $E_p(\sigma)$,
it is easy to see that
$$\eqalign{E''_p(\sigma)E_p(\sigma)-E'_p(\sigma)^2
& = {4\over \pi^2}
\int_0^{\pi} D_p(\theta)^\sigma \log^2D_p(\theta) \sin^2\theta \d\theta
\int_0^{\pi} D_p(\theta)^\sigma \sin^2\theta \d\theta
\cr
& \quad
- \bigg({2\over \pi}
\int_0^{\pi} D_p(\theta)^\sigma 
\log D_p(\theta) \sin^2\theta \d\theta\bigg)^2
\cr\noalign{\smallskip}
& = {4\over \pi^2}
\int_0^{\pi} \int_0^{\pi}
D_p(\theta_1)^\sigma D_p(\theta_2)^\sigma 
\big(\log^2D_p(\theta_1) - \log D_p(\theta_1) \log D_p(\theta_2)\big)
\times
\cr
& \qquad\times
\sin^2\theta_1 \sin^2\theta_2 \d\theta_1 \d\theta_2.
\cr}$$
In view of the symmetry in $\theta_1$ and $\theta_2$,
the same formula holds 
if we exchange the roles of $\theta_1$ and $\theta_2$.
Thus it follows that
$$
E''_p(\sigma)E_p(\sigma)-E'_p(\sigma)^2
= {2\over \pi^2}
\int_0^{\pi} \int_0^{\pi}
D_p(\theta_1)^\sigma D_p(\theta_2)^\sigma
\log^2\bigg({D_p(\theta_1)\over D_p(\theta_2)}\bigg) 
\sin^2\theta_1 \sin^2\theta_2 \d\theta_1 \d\theta_2.
$$
This proves the first assertion
and the second follows immediately.
\hfill
$\square$

\proclaim Lemma 2.3.
There is an absolute constant $c\ge 2$ such that
for $t\ge 4\log c$ and $y\ge ce^t$, 
the equation $\phi_1(\sigma, y)=2(\log t+\gamma)$
has a unique positive solution in $\sigma$.
Denoting by $\kappa(t, y)$ this solution, 
we have 
$\kappa(t, y)\asymp e^t$
uniformly for $t\ge 4\log c$ and $y\ge ce^t$.

\noindent{\sl Proof}.
According to Lemma 4.3 below with the choice of $J=1$, we have 
$$\phi_1(\sigma, y)=2(\log_2\sigma+\gamma) + O(1/\log\sigma)$$
for $y\ge \sigma\ge 2$.
Thus
$$\eqalign{\phi(ce^t, y)
& = 2\log(t+\log c) + 2\gamma + O\bigg({1\over t+\log c}\bigg)
\cr
& > 2\log t + 2\gamma
\cr}$$
and
$$\eqalign{
\phi(c^{-1}e^t, y)
& = 2\log(t-\log c) + 2\gamma + O\bigg({1\over t-\log c}\bigg)
\cr
& < 2\log t + 2\gamma,
\cr}$$
provided that $c$ is a large constant and $t\ge 4\log c$.
On the other hand,
in view of Lemma~2.2, 
we know that for any $y\ge 2$,
$\phi_1(\sigma, y)$ is an increasing function of $\sigma$ 
in $(0, \infty)$.
Hence the equation $\phi_1(\sigma, y)=2(\log t+\gamma)$
has a unique positive solution $\kappa(t, y)$ and
$c^{-1}e^t\le \kappa(t, y)\le ce^t$
for $t\ge 4\log c$ and $y\ge ce^t$.
This completes the proof.
\hfill
$\square$

\goodbreak
\vskip 8mm

\centerline{\pptitre \S\ 3. Preliminary lemmas}

\bigskip

This section is devoted to establish some preliminary lemmas,
which will be useful later.

\proclaim Lemma 3.1.
Let $j\ge 0$ be a fixed real number.
Then we have
$$\int_0^{\pi} e^{2u\cos\theta} 
(1-\cos\theta)^j \sin^2\theta \d\theta
\asymp_j e^{2u}u^{-(j+3/2)}
\qquad(u\ge 1).
\leqno(3.1)$$
The implied constant depends on $j$ only.

\noindent{\sl Proof}.
First we write
$$\eqalign{
\int_0^{\pi} e^{2u\cos\theta} 
(1-\cos\theta)^j \sin^2\theta \d\theta
& = \int_0^{\pi/2} 
\big(e^{2u\cos\theta} (1-\cos\theta)^j 
+ e^{-2u\cos\theta} (1+\cos\theta)^j \big)
\sin^2\theta \d\theta
\cr
& = \int_0^{1} 
\big(e^{2ut} (1-t)^j 
+ e^{-2ut} (1+t)^j\big) 
(1-t^2)^{1/2} \d t
\cr
& \asymp \int_0^{1} e^{2ut} (1-t)^{j+1/2} \d t
+ \int_0^{1} e^{-2ut} (1-t)^{1/2} \d t.
\cr}$$
By the change of variables $u(1-t)=v$, it follows that
$$\eqalign{
\int_0^{1} e^{2ut} (1-t)^{j+1/2} \d t
& = e^{2u}u^{-(j+3/2)} 
\int_0^{u} e^{-2v} v^{j+1/2} \d v
\cr
& \asymp e^{2u}u^{-(j+3/2)},
\cr
\int_0^{1} e^{-2ut} (1-t)^{1/2} \d t
& \le \int_0^{1} e^{-2ut} \d t
\ll u^{-1}.
\cr}$$ 
We obtain the desired result 
by insertion of these estimates into the preceeding relation.
\hfill
$\square$

\proclaim Lemma 3.2.
Let $j\ge 0$ be an integer and
$$E_{p,j}(\sigma)
:= {2\over \pi}
\int_0^{\pi} D_p(\theta)^\sigma (1-\cos\theta)^j \sin^2\theta \d\theta.
\leqno(3.2)$$
{\rm (}In particular $E_{p,0}(\sigma)=E_p(\sigma).)$
Then we have
$$E_{p,j}(\sigma)
= {2^{j+3}\over \pi} \int_0^{1}
\bigg[\bigg(1-{1\over p}\bigg)^2+{4u\over p}\bigg]^{-\sigma}
u^{j+1/2} (1-u)^{1/2} \d u$$
and the estimate
$$E_{p,j}(\sigma)/E_p(\sigma)\ll (p/\sigma)^j
\leqno(3.3)$$
holds uniformly for all primes $p$ and $\sigma>0$.
Further if $p\ge \sigma\ge 0$, we have
$$E_p(\sigma)
\asymp 1.
\leqno(3.4)$$
The implied constant in (3.3) depends on $j$ only
and the one in (3.4) is absolute.

\noindent{\sl Proof}.
By the change of variables $u=\sin^2(\theta/2)$, 
a simple computation shows that the first assertion is true.
Obviously (3.3) holds for $j=0$.

Now assume that it is true for $j$.
An integration by parts leads to
$$\eqalign{
E_p(\sigma)
& \gg_j \bigg({\sigma\over p}\bigg)^j 
\int_0^{1}
\bigg[\bigg(1-{1\over p}\bigg)^2+{4u\over p}\bigg]^{-\sigma}
u^{j+1/2} (1-u)^{1/2} \d u
\cr
& \gg_j \bigg({\sigma\over p}\bigg)^j 
\int_0^{1}
\bigg\{
\bigg[\bigg(1-{1\over p}\bigg)^2+{4u\over p}\bigg]^{-1}
{4\sigma\over p}
+ {1\over 2(1-u)}
\bigg\}
\times
\cr
& \qquad\times
\bigg[\bigg(1-{1\over p}\bigg)^2+{4u\over p}\bigg]^{-\sigma}
u^{j+1+1/2} (1-u)^{1/2} \d u.
\cr}$$
On the other hand, we have
$$0<u<1\;\;\Rightarrow\;\;
\bigg[\bigg(1-{1\over p}\bigg)^2+{4u\over p}\bigg]^{-1}
{4\sigma\over p}
+ {1\over 2(1-u)}
\ge \bigg(1+{1\over p}\bigg)^{-2}
{4\sigma\over p}
\ge {16\sigma\over 9p}.$$
Inserting it into the preceeding estimate, 
we see that 
$$\eqalign{
E_p(\sigma)
& \gg_j \bigg({\sigma\over p}\bigg)^{j+1} 
\int_0^{1}
\bigg[\bigg(1-{1\over p}\bigg)^2+{4u\over p}\bigg]^{-\sigma}
u^{j+1+1/2} (1-u)^{1/2} \d u
\cr
& \asymp_j \bigg({\sigma\over p}\bigg)^{j+1} E_{p, j+1}(\sigma).
\cr}$$
Thus (3.3) holds also for $j+1$.

Since
$(1+1/p)^{-2}
\le D_p(\theta)
\le (1-1/p)^{-2}$
for all primes $p$ and any $\theta\in \R$,
we have
$D_p(\theta)^\sigma\asymp 1$
uniformly for $p\ge \sigma\ge 0$ and $\theta\in \R$.
This implies (3.4).
\hfill
$\square$

\medskip

Introduce the function
$$g(u) 
:= \log\bigg({2\over \pi}
\int_0^{\pi} e^{2u\cos\theta} \sin^2\theta \d\theta
\bigg)
\qquad
(u\ge 0)
\leqno(3.5)$$
and let $h(u)$ be defined as in (1.22).
Clearly we have
$$\leqalignno{
h(u) 
& = \cases{
g(u)       & if $\;\;0\le u<1$,
\cr\noalign{\smallskip}
g(u) - 2u  & if $\;\;u\ge 1$,
\cr}
& (3.6)
\cr\noalign{\smallskip}
h'(u) 
& = \cases{
g'(u)      & if $\;\;0\le u<1$,
\cr\noalign{\smallskip}
g'(u) - 2  & if $\;\;u\ge 1$,
\cr}
& (3.7)
\cr\noalign{\medskip}
h''(u)
& = g''(u)
\quad(u\ge 0, \, u\not=1).
& (3.8)
\cr}$$

\proclaim Lemma 3.3.
We have
$$\leqalignno{
h(u)
& \asymp \cases{
u^2        & if $\;\;0\le u<1$,
\cr\noalign{\smallskip}
\log(2u)   & if $\;\;u\ge 1$,
\cr}
& (3.9)
\cr\noalign{\smallskip}
h'(u)
& \asymp \cases{
u        & if $\;\;0\le u<1$,
\cr\noalign{\smallskip}
u^{-1}   & if $\;\;u\ge 1$,
\cr}
& (3.10)
\cr\noalign{\smallskip}
h''(u)
& \asymp \cases{
1        & if $\;\;0\le u<1$,
\cr\noalign{\smallskip}
u^{-2}   & if $\;\;u\ge 1$,
\cr}
& (3.11)
\cr\noalign{\smallskip}
h'''(u)
& \asymp \cases{
u        & if $\;\;0\le u<1$,
\cr\noalign{\smallskip}
u^{-3}   & if $\;\;u\ge 1$.
\cr}
& (3.12)
\cr}$$

\noindent{\sl Proof}.
When $0\le u<1$, we have
$$e^{2u\cos\theta}
= \sum_{n=0}^\infty {(u\cos\theta)^n\over n!}.$$
From this we deduce that
$$\leqalignno{h(u)
& =\log\bigg(
{2\over \pi} \sum_{n=0}^\infty  {u^n\over n!}
\int_0^{\pi} (\cos\theta)^n \sin^2\theta \d\theta
\bigg)
& (3.13)
\cr
& =\log\bigg(1 +
\sum_{\ell=1}^\infty  
{2\cdot (2\ell-1)!!\over (2\ell)! (2\ell+2)!!} u^{2\ell}
\bigg),
\cr}$$
where we have used the following facts:
$$\int_0^{\pi} (\cos\theta)^{2\ell+1} \sin^2\theta \d\theta
=0$$
and
$${2\over \pi} 
\int_0^{\pi} (\cos\theta)^{2\ell} \sin^2\theta \d\theta
= \cases{
1                               & if $\;\ell=0$,
\cr\noalign{\medskip}
\displaystyle
2{(2\ell-1)!!\over (2\ell+2)!!} & if $\;\ell\ge 1$
\cr}$$
and where $n!!$ denotes the product of all positive integer from 1 to $n$ %
having same parity than $n$.
Now we easily deduce, from (3.13), the desired results (3.9)--(3.12) 
in the case of $0\le u<1$.

The estimates of (3.9)--(3.12) for $u>1$ are simple consequences of (3.1),
by noticing the following relations
$$\eqalign{h'(u)
& = - 2\,{\displaystyle
\int_0^{\pi} e^{2u\cos\theta} (1-\cos\theta) \sin^2\theta \d\theta
\over 
\displaystyle
\int_0^{\pi} e^{2u\cos\theta} \sin^2\theta \d\theta},
\cr\noalign{\smallskip}
h''(u)
& = 4\,{\displaystyle
\int_0^{\pi} e^{2u\cos\theta} (1-\cos\theta)^2 \sin^2\theta \d\theta
\over 
\displaystyle
\int_0^{\pi} e^{2u\cos\theta} \sin^2\theta \d\theta}
- 4 \left(
{\displaystyle
\int_0^{\pi} e^{2u\cos\theta} (1-\cos\theta) \sin^2\theta \d\theta
\over 
\displaystyle
\int_0^{\pi} e^{2u\cos\theta} \sin^2\theta \d\theta}
\right)^2.
\cr}$$
This completes the proof.
\hfill
$\square$

\vskip 8mm

\centerline{\pptitre \S\ 4. Estimates of $\phi_n(\sigma, y)$}

\bigskip

The aim of this section is to prove some estimates of $\phi_n(\sigma, y)$
for $n=0, 1, 2, 3, 4$.

\proclaim Lemma 4.1.
For any fixed integer $J\ge 1$, we have
$$\phi_0(\sigma, y)
= \sigma 
\bigg\{
2\log_2\sigma + 2\gamma
+ \sum_{j=1}^J {b_{j, 0}\over (\log\sigma)^j}
+ O_J\big(R_J(\sigma, y)\big)\bigg\}
\leqno(4.1)$$
uniformly for $y\ge \sigma\ge 3$,
where $R_J(\sigma, y)$ is defined as in (1.20) and
$$b_{j, 0} := \int_0^\infty {h(u)\over u^2} (\log u)^{j-1} \d u.
\leqno(4.2)$$

\noindent{\sl Proof}.
By the definition (2.1) of $D_p(\theta)$ and the one of $E_p(\sigma)$,
it is easy to see that
for $p\ge \sigma^{1/2}$, we have
$$\leqalignno{D_p(\theta)^\sigma
& = e^{2(\sigma/p)\cos\theta}
\bigg\{1+O\bigg({\sigma\over p^2}\bigg)\bigg\},
& (4.3)
\cr\noalign{\medskip}
E_p(\sigma) 
& = \bigg\{1+O\bigg({\sigma\over p^2}\bigg)\bigg\}
{2\over \pi} \int_0^{\pi} e^{2(\sigma/p)\cos\theta}
\sin^2\theta \d\theta.
& (4.4)
\cr}$$
From these, we deduce that
$$\sum_{\sigma^{1/2}<p\le y} \log E_p(\sigma)
= \sum_{\sigma^{1/2}<p\le y} g(\sigma/p)
+ O(\sigma^{1/2}/\log\sigma)
\leqno(4.5)$$
where $g(u)$ is defined as in (3.5).

In order to treat the sum over $p\le \sigma$, we write
$$E_p(\sigma)
= (1 - 1/p)^{-2\sigma}
E_p^*(\sigma),$$
where
$$E_p^*(\sigma)
:= {2\over \pi} \int_0^{\pi} 
\bigg\{1 + {2(1-\cos\theta)\over p}\bigg(1-{1\over p}\bigg)^{-2}
\bigg\}^{-\sigma}
\sin^2\theta \d\theta.$$
By using the change of variables $u = \sin^2(\theta/2)$, we have
$$\eqalign{E_p^*(\sigma)
& = {8\over \pi} \int_0^{\pi} 
\bigg\{1 + {4\over p}\bigg(1-{1\over p}\bigg)^{-2}\sin^2(\theta/2)
\bigg\}^{-\sigma}
\sin^2(\theta/2)\cos^2(\theta/2) \d\theta
\cr
& \ge {8\over \pi}
\int_0^{p/2\sigma} 
\bigg\{1 + {4\over p}\bigg(1-{1\over p}\bigg)^{-2}u
\bigg\}^{-\sigma}
\sqrt{u (1-u)} \d u
\cr
& \ge {8\over \pi}
\bigg(1 + {8\over \sigma}\bigg)^{-\sigma}
\int_0^{p/2\sigma} \sqrt{u (1-u)} \d u
\cr
& \ge C \bigg({p\over \sigma}\bigg)^{3/2},
\cr}$$
where $C>0$ is a constant.
On the other hand, we have trivially $E_p^*(\sigma)\le 1$
for all $p$ and $\sigma>0$.
Thus $|\log E_p^*(\sigma)|\ll \log(\sigma/p)$ for $p\le \sigma^{1/2}$
and
$$\sum_{p\le \sigma^{1/2}} \big|\log E_p^*(\sigma)\big|
\ll \sum_{p\le \sigma^{1/2}} \log(\sigma/p)
\ll \sigma^{1/2}.
\leqno(4.6)$$
Combining (4.5) and (4.6), we can write
$$\sum_{p\le y} \log E_p(\sigma)
= 2\sigma \sum_{p\le \sigma^{1/2}} \log(1-1/p)^{-1}
+ \sum_{\sigma^{1/2}<p\le y} g(\sigma/p)
+ O(\sigma^{1/2}).$$

In view of (3.6) and the following estimate
$$\sum_{\sigma^{1/2}<p\le \sigma} 
\big(2\sigma \log(1-1/p)^{-1} - 2\sigma/p\big)
\ll \sum_{\sigma^{1/2}<p\le \sigma} \sigma/p^2
\ll \sigma^{1/2}/\log\sigma,$$
the preceeding estimate can be written as
$$\sum_{p\le y} \log E_p(\sigma)
= 2\sigma \sum_{p\le \sigma} \log(1-1/p)^{-1}
+ \sum_{\sigma^{1/2}<p\le y} h(\sigma/p)
+ O(\sigma^{1/2}).
\leqno(4.7)$$

By using the prime number theorem in the form
$$\pi(t) 
:= \sum_{p\le t} 1
= \int_2^t {\d v\over \log v}
+ O\Big(te^{-8\sqrt{\log t}}\Big),
\leqno(4.8)$$
it follows that 
$$\sum_{\sigma^{1/2}<p\le y} h\bigg({\sigma\over p}\bigg)
= \int_{\sigma^{1/2}}^y 
{h(\sigma/t)\over \log t} \d t
+ O(R_0),
\leqno(4.9)$$
where
$$\eqalign{R_0
& := h\bigg({\sigma\over y}\bigg)ye^{-8\sqrt{\log y}}
+ h\big(\sigma^{1/2}\big)\sigma^{1/2}e^{-4\sqrt{\log \sigma}}
+ \int_{\sigma^{1/2}}^y 
(\sigma/t)|h'(\sigma/t)|e^{-8\sqrt{\log t}} \d t
\cr
& \,\ll {\sigma^2\over y} e^{-8\sqrt{\log y}}
+ \sigma^{1/2}e^{-2\sqrt{\log\sigma}}
+ \int_{\sigma^{1/2}}^{\sigma} 
e^{-2\sqrt{\log t}}\d t
+ \sigma^2 \int_{\sigma}^y 
{e^{-8\sqrt{\log t}}\over t^2} \d t
\cr
& \,\ll \sigma e^{-\sqrt{\log\sigma}}
\cr}$$
by use of Lemma 3.3.

In order to evaluate the integral of (4.9), 
we use the change of variables $u = \sigma/t$ to write
$$\eqalign{\int_{\sigma^{1/2}}^y {h(\sigma/t)\over \log t} \d t
& = \sigma \int_{\sigma/y}^{\sigma^{1/2}} 
{h(u)\over u^2\log(\sigma/u)} \d u
\cr
& = \sigma \int_{\sigma^{-1/2}}^{\sigma^{1/2}} 
{h(u)\over u^2\log(\sigma/u)} \d u
+ O\big(R_0'\big)
\cr}$$
where
$$\eqalign{R_0'
& := \sigma \int_0^{\sigma/y} 
{|h(u)|\over u^2\log(\sigma/u)} \d u
+ \sigma \int_0^{\sigma^{-1/2}} 
{|h(u)|\over u^2\log(\sigma/u)} \d u
\cr
& \,\ll {\sigma^2\over y\log y}
+ {\sigma^{1/2}\over \log\sigma}.
\cr}$$
On the other hand, we have
$$\eqalign{\int_{\sigma^{-1/2}}^{\sigma^{1/2}} 
{h(u)\over u^2\log(\sigma/u)} \d u
& = {1\over \log\sigma}
\int_{\sigma^{-1/2}}^{\sigma^{1/2}} 
{h(u)\over u^2(1-(\log u)/\log\sigma)} \d u
\cr
& = \sum_{j=1}^J {1\over (\log\sigma)^j}
\int_{\sigma^{-1/2}}^{\sigma^{1/2}} 
{h(u)\over u^2} (\log u)^{j-1} \d u
+ O\bigg({1\over (\log\sigma)^{J+1}}\bigg).
\cr}$$
Extending the interval of integration $[\sigma^{-1/2}, \, \sigma^{1/2}]$ 
to $(0, \infty)$ and bounding the contributions of 
$(0, \sigma^{-1/2}]$ and $[\sigma^{1/2}, \infty)$ 
by using (3.9) of Lemma 3.3,
we have
$$\int_{\sigma^{-1/2}}^{\sigma^{1/2}} {h(u)\over u^2} (\log u)^{j-1} \d u
= b_{j, 0} 
+ O\bigg(
{{(\log\sigma)^j}\over{ \sigma^{1/2}}}
\bigg).
$$
Combining these estimates, we find that
$$\sum_{\sigma^{1/2}<p\le y} h\bigg({\sigma\over p}\bigg) 
= \sigma \bigg\{
\sum_{j=1}^J {b_{j, 0}\over (\log\sigma)^j}
+ O_J\big(R_J(\sigma, y)\big)
\bigg\}.
\leqno(4.10)$$
Now the desired result follows from (4.7), (4.10) and
 the prime number theorem in the form
$$\sum_{p\le \sigma} \log(1-1/p)^{-1}
= \log_2\sigma + \gamma 
+ O\Big(e^{-2\sqrt{\log\sigma}}\Big).
\leqno(4.11)$$
This completes the proof.
\hfill
$\square$

\medskip

{\bf Remark 3}.
In view of (1.3), we can write (4.1) as
$$\phi_0(\sigma, y)
= \sigma 
\bigg\{
\log(B_1^+\log\sigma)^{A_1^+}
+ \sum_{j=1}^J {b_{j, 0}\over (\log\sigma)^j}
+ O_J\big(R_J(\sigma, y)\big)\bigg\}$$
uniformly for $y\ge \sigma\ge 3$.
In the case $\sigma<0$,
a similar asymptotic formula 
(with $A_1^-, B_1^-$ and corresponding $b_{j, 0}^-$ 
in place of $A_1^+, B_1^+$ and $b_{j, 0}$)
can be established uniformly for $y\ge -\sigma\ge 3$.
As indicated in the introduction,
Lemma 4.1 can be easily generalised to the general case $m\ge 1$.
Thus we give an improvement and generalisation 
of Corollaries A and C of [15], 
of Theorem B of [5],
and an improvement of Theorem 1.12 of [1].
It is worthy to indicate that 
our method seems to be simpler and more natural.

\smallskip

\proclaim Lemma 4.2.
We have
$${E_p'(\sigma)\over E_p(\sigma)}
= \cases{
\displaystyle
\log D_p(0) 
+ O\bigg({1\over \sigma}\bigg)
& for all $p$ and $\sigma>0$, 
\cr\noalign{\bigskip}
\displaystyle
{1\over 2}
g'\bigg({\sigma\over p}\bigg) \log D_p(0)
+ O\bigg({1\over p^2} + {\sigma\over p^3}\bigg)
& if $\;\;p\ge \sigma^{1/2}$,
\cr}
\leqno(4.12)$$
where $g(u)$ is defined as in (3.5).

\noindent{\sl Proof}.
First we write
$$\leqalignno{E_p'(\sigma)
& = {2\over \pi} \int_0^{\pi} 
D_p(\theta)^{\sigma} \log D_p(\theta)
\sin^2\theta \d\theta
& (4.13)
\cr
& = E_p(\sigma) \log D_p(0) 
+ R',
\cr}$$
where
$$R'
:= {2\over \pi} \int_0^{\pi} 
D_p(\theta)^{\sigma}
\log\bigg({D_p(\theta)\over D_p(0)}\bigg)
\sin^2\theta \d\theta.
\leqno(4.14)$$
Since
$$\bigg|\log\bigg({D_p(\theta)
\over D_p(0)}\bigg)\bigg|
= \bigg|-\log\bigg(1+{2p(1-\cos\theta)\over (p-1)^2}\bigg)\bigg|
\le {2p(1-\cos\theta)\over (p-1)^2}
\le {8(1-\cos\theta)\over p},$$
it follows from (3.3) of Lemma 3.2 with $j=1$ that
$${R'\over E_p(\sigma)}
\ll {E_{p,1}(\sigma)\over p E_p(\sigma)}
\ll {1\over \sigma}$$
for all $p$ and $\sigma>0$.
This implies, via (4.13), the first estimate of (4.12).

\smallskip

We have
$$\eqalign{\log D_p(\theta)
& = (\cos\theta) (2/p) + O\big(1/p^2\big)
\cr
& = (\cos\theta) \log D_p(0) + O\big(1/p^2\big).
\cr}$$
Inserting it and (4.3) into the first relation of (4.13)
and in view of (4.4),
we can write, for $p\ge \sigma^{1/2}$,
$$\leqalignno{E_p'(\sigma) 
& = \bigg\{1+O\bigg({\sigma\over p^2}\bigg)\bigg\}
{2\over \pi} \int_0^{\pi} e^{2(\sigma/p)\cos\theta}
\bigg\{(\cos\theta) \log D_p(0) 
+ O\bigg({1\over p^2}\bigg)\bigg\}
\sin^2\theta \d\theta
\cr
& = \bigg\{1+O\bigg({\sigma\over p^2}\bigg)\bigg\} 
{2\over \pi}\int_0^{\pi} e^{2(\sigma/p)\cos\theta}
(\cos\theta) \sin^2\theta \d\theta
\log D_p(0) 
+ O\bigg({E_p(\sigma)\over p^2}\bigg).
\cr}$$
From this and (4.4), we deduce
$${E_p'(\sigma)\over E_p(\sigma)}
= \bigg\{1+O\bigg({\sigma\over p^2}\bigg)\bigg\}
{1\over 2}
g'\bigg({\sigma\over p}\bigg) \log D_p(0)
+ O\bigg({1\over p^2}\bigg),
$$
which implies the second estimate of (4.12).
This completes the proof.
\hfill
$\square$

\smallskip

\proclaim Lemma 4.3.
Let $J\ge 1$ be a fixed integer.
Then we have
$$\phi_1(\sigma, y)
= 2\log_2\sigma + 2\gamma
+ \sum_{j=1}^J {b_{j, 1}\over (\log\sigma)^j}
+ O_J\big(R_J(\sigma, y)\big)$$
uniformly for $y\ge \sigma\ge 3$, 
where the constant $b_{j, 1}$ is given by
$$b_{j, 1} := \int_0^\infty {h'(u)\over u} (\log u)^{j-1} \d u
\leqno(4.15)$$
and $R_J(\sigma, y)$ is defined as in (1.20).

\noindent{\sl Proof}.
We have
$$\phi_1(\sigma, y)
= \sum_{p\le y} {E_p'(\sigma)/E_p(\sigma)}.$$
Using the first relation of (4.12) 
for $p\le \sigma^{2/3}$ 
and the second for $\sigma^{2/3}<p\le y$, 
we obtain
$$\phi_1(\sigma, y)
= \sum_{p\le \sigma^{2/3}} \log D_p(0)
+ {1\over 2}\sum_{\sigma^{2/3}<p\le y}
g'\bigg({\sigma\over p}\bigg) \log D_p(0)
+ O\bigg({1\over \sigma^{1/3}}\bigg).$$
In view of (3.7), 
the preceeding formula can be written as
$$\phi_1(\sigma, y)
= \sum_{p\le \sigma} \log D_p(0)
+ \sum_{\sigma^{2/3}<p\le y} 
h'\bigg({\sigma\over p}\bigg) \log\bigg(1-{1\over p}\bigg)^{-1} 
+ O\bigg({1\over \sigma^{1/3}}\bigg).
\leqno(4.16)$$

Similarly to (4.10), we can prove that
$$\sum_{\sigma^{2/3}<p\le y} 
h'\bigg({\sigma\over p}\bigg) \log\bigg(1-{1\over p}\bigg)^{-1}
= \sum_{j=1}^J {b_{j, 1}\over (\log\sigma)^j}
+ O_J\big(R_J(\sigma, y)\big),
\leqno(4.17)$$
using (3.10), (3.11) and (4.11) instead of (3.9), (3.10) and (4.8).
Now the desired result follows from (4.16), (4.10) and (4.17).
\hfill
$\square$

\smallskip

\proclaim Lemma 4.4.
We have
$${E_p''(\sigma)E_p(\sigma)-E_p'(\sigma)^2\over E_p(\sigma)^2}
= \cases{\displaystyle
O\bigg({1\over \sigma^2}\bigg)                
& if $\;\;p\le \sigma^{1/2}$,
\cr\noalign{\medskip}
\displaystyle
{1\over p^2}g''\bigg({\sigma\over p}\bigg) 
+ O\bigg(\min\bigg\{{1\over \sigma^{2}p}, \, 
{1\over \sigma p^2}\bigg\}\bigg)  
& if $\;\;p>\sigma^{1/2}$,
\cr}
\leqno(4.18)$$
where $g(u)$ is defined as in (3.5). 

\noindent{\sl Proof}.
First we write
$$\leqalignno{E''_p(\sigma)
& = {2\over \pi} \int_0^{\pi} 
D_p(\theta)^{\sigma}
\log^2D_p(\theta)
\sin^2\theta \d\theta
& (4.19)
\cr
& = E_p(\sigma) \log^2D_p(0)
+ R'', 
\cr}$$
where
$$R''
:= {2\over \pi} \int_0^{\pi} 
D_p(\theta)^{\sigma}
\Big(
\log^2D_p(\theta)
- \log^2D_p(0)
\Big)
\sin^2\theta \d\theta.$$
Using (4.13) and (4.19), we can deduce 
$${E_p''(\sigma)E_p(\sigma)-E_p'(\sigma)^2\over E_p(\sigma)^2}
= {R'' - 2R'\log D_p(0)\over  E_p(\sigma)}
- \bigg({R'\over E_p(\sigma)}\bigg)^2,
\leqno(4.20)$$
where $R'$ is defined as in (4.14).

From the definitions of $R'$ and $R''$,
a simple calculation shows that
$$R'' - 2 R' \log D_p(0)
= {2\over \pi} \int_0^{\pi} 
D_p(\theta)^{\sigma}
\log^2
\bigg({D_p(\theta)\over D_p(0)}\bigg)
\sin^2\theta \d\theta.$$
Since
$$\log^2\bigg({D_p(\theta)\over D_p(0)}\bigg)
= \log^2\bigg(1+{2p(1-\cos\theta)\over (p-1)^2}\bigg)
= {4(1-\cos\theta)^2\over p^2}
+ O\bigg({(1-\cos\theta)^2\over p^3}\bigg),$$
we have
$$R'' - 2 R' \log D_p(0)
= {4\over p^2} E_{p,2}(\sigma)
+ O\bigg({E_{p,2}(\sigma)\over p^3}\bigg),$$
where $E_{p,j}(\sigma)$ is defined as in (3.2).
By using (3.3) with the choice of $j=2$
and the trivial estimate $E_{p,2}(\sigma)\le 4E_p(\sigma)$, 
we deduce 
$${R'' - 2 R' \log D_p(0)\over E_p(\sigma)}
= {4\over p^2} {E_{p,2}(\sigma)\over  E_p(\sigma)}
+ O\bigg(\min\bigg\{{1\over \sigma^2p}, \; 
{1\over p^3}\bigg\}\bigg).
\leqno(4.21)$$

Similarly we have
$$\log\bigg({D_p(\theta)\over D_p(0)}\bigg)
= -\log\bigg(1+{2p(1-\cos\theta)\over (p-1)^2}\bigg)
= -{2(1-\cos\theta)\over p}
+ O\bigg({(1-\cos\theta)\over p^2}\bigg),$$
and therefore
$$R'
= - {2\over p} E_{p,1}(\sigma)
+ O\bigg({E_{p,1}(\sigma)\over p^2}\bigg).$$
Now (3.3) with $j=1$ and 
the trivial estimate $E_{p, 1}(\sigma)\le 2E_p(\sigma)$ imply 
$$\leqalignno{\bigg({R'\over E_p(\sigma)}\bigg)^2
& = {4\over p^2} \bigg({E_{p,1}(\sigma)\over  E_p(\sigma)}\bigg)^2
+ O\bigg({E_{p,1}(\sigma)^2\over p^3E_p(\sigma)^2}\bigg)
& (4.22)
\cr
& = {4\over p^2} \bigg({E_{p,1}(\sigma)\over  E_p(\sigma)}\bigg)^2
+ O\bigg(\min\bigg\{{1\over \sigma^2p}, \; 
{1\over p^3}\bigg\}\bigg).
\cr}$$
Inserting (4.21) and (4.22) into (4.20) 
and in view of (4.14),
we deduce
$${E_p''(\sigma)E_p(\sigma)-E_p'(\sigma)^2\over E_p(\sigma)^2}
= {4\over p^2}h_p(\sigma) 
+ O\bigg(\min\bigg\{{1\over \sigma^2 p}, \; 
{1\over p^3}\bigg\}\bigg)
\leqno(4.23)$$
for all $p$ and $\sigma>0$,
where
$$h_p(\sigma)
:= {E_{p,2}(\sigma)\over  E_p(\sigma)}
- \bigg({E_{p,1}(\sigma)\over  E_p(\sigma)}\bigg)^2.$$

When $p\le \sigma^{1/2}$, the inequality (3.3) of Lemma 3.2 implies that
$h_p(\sigma)\ll (p/\sigma)^2$.
From this and (4.23) we deduce the first estimate of (4.18).

If $p\ge \sigma^{1/2}$, 
we can use (4.3), (3.11) and (3.8) to write
$$\eqalign{4h_p(\sigma)
& = g''\bigg({\sigma\over p}\bigg)
\bigg\{1+O\bigg({\sigma\over p^2}\bigg)\bigg\}
\cr\noalign{\smallskip}
& = g''\bigg({\sigma\over p}\bigg)
+ O\bigg(\min\bigg\{{\sigma \over p^2}, 
\, {1\over \sigma}\bigg\}\bigg).
\cr}$$
Inserting it into (4.23) and in view of Lemma 3.1, 
we get, for $p\ge \sigma^{1/2}$, 
$$\leqalignno{\qquad
{E_p''(\sigma)E_p(\sigma)-E_p'(\sigma)^2\over E_p(\sigma)^2}
& = {1\over p^2}g''\bigg({\sigma\over p}\bigg)
+ O\bigg(\min\bigg\{{1\over \sigma^2 p},\,
{1\over \sigma p^2}\bigg\}\bigg).
\cr}$$
This completes the proof.
\hfill
$\square$

\proclaim Lemma 4.5.
Let $J\ge 1$ be a fixed integer.
Then we have
$$\phi_2(\sigma, y)
= {1\over \sigma}
\bigg\{\sum_{j=1}^J {b_{j, 2}\over (\log\sigma)^j}
+ O_J\big(R_J(\sigma, y)\big)\bigg\}$$
uniformly for $y\ge \sigma\ge 2$,
where
$$b_{j,2}
:= \int_{0}^{\infty} h''(u) (\log u)^{j-1} \d u.$$
In particular $b_{1, 2}=2$.

\noindent{\sl Proof}.
From Lemma 4.4 and (3.8), we deduce easily that
$$\eqalign{\phi_2(\sigma, y)
& = \sum_{p\le y}
{E_p''(\sigma)E_p(\sigma)-E_p'(\sigma)^2\over E_p(\sigma)^2}
\cr
& = \sum_{\sigma^{1/2}<p\le y} {g''(\sigma/p)\over p^2} 
+ O\bigg({1\over \sigma^{3/2}\log\sigma}\bigg)
\cr
& = \sum_{\sigma^{1/2}<p\le y} {h''(\sigma/p)\over p^2} 
+ O\bigg({1\over \sigma^{3/2}\log\sigma}\bigg).
\cr}$$
Similarly to (4.10), we can prove that
$$\sum_{\sigma^{1/2}<p\le y} {h''(\sigma/p)\over p^2} 
= {1\over \sigma}
\bigg\{\sum_{j=1}^J {b_{j, 2}\over (\log\sigma)^j}
+ O_J\big(R_J(\sigma, y)\big)\bigg\},$$
by using (3.11), (3.12) and (4.8).
Now the desired result follows from the preceeding two estimates.

Finally
$$\eqalign{b_{1, 2}
& = \int_{0}^{1} h''(u) \d u + \int_{1}^{\infty} h''(u) \d u
\cr
& = h'(1-) - h'(1+)
= h'(1-) - \big(h'(1-)-2\big)
= 2.
\cr}$$
This completes the proof.
\hfill
$\square$

\smallskip

Similarly (even more easily, 
since we only need an upper bound instead of an asymptotic formula), 
we can prove the following result.

\proclaim Lemma 4.6.
We have
$$\phi_n(\sigma, y)
\ll 1/(\sigma^{n-1}\log\sigma)
\qquad(n=3, 4)
\leqno(4.24)$$
uniformly for $y\ge \sigma\ge 3$.

\vskip 8mm

\centerline{\pptitre \S\ 5. Estimate of $|E(\kappa+i\tau, y)|$}

\bigskip

\proclaim Lemma 5.1.
For any $\delta\in (0, {1\over 4})$,
there are two absolute positive constants $c_1, c_2$ 
and a positive constant $c_3=c_3(\delta)$ such that
for all $y\ge \sigma\ge 3$ we have
$$\bigg|{E(\sigma+i\tau, y)\over E(\sigma, y)}\bigg|
\le \cases{
1 
& if $\;\;|\tau|\le c_1\sigma^{1/2}\log\sigma\;$
or $\;|\tau|\ge y^{1/\delta}$,
\cr\noalign{\smallskip}
e^{-c_2\tau^2/[\sigma(\log\sigma)^2]} 
& if $\;\;c_1\sigma^{1/2}\log\sigma\le |\tau|\le \sigma$,
\cr\noalign{\smallskip}
e^{-c_3|\tau|^\delta} 
& if $\;\;\sigma\le |\tau|\le y^{1/\delta}$.
\cr}
\leqno(5.1)$$

\noindent{\sl Proof}.
First we write
$$\eqalign{E_p(s)
& = {2\over \pi} 
\int_0^\pi \big(D_p(\theta)^{-1}\big)^{-s} 
\sin^2\theta \d \theta
\cr
& = {2\over \pi} 
\int_0^\pi {\sin^2\theta\over (1-s)(D_p(\theta)^{-1})'}
\d \big(D_p(\theta)^{-1}\big)^{1-s}. 
\cr}$$
Since $(D_p(\theta)^{-1})'=2p^{-1}\sin\theta$,
after a simplification and an integration by parts
it follows that
$$\eqalign{E_p(s)
& = {p\over \pi (s-1)}
\int_0^\pi D_p(\theta)^{s-1} \cos\theta \d \theta
\cr
& = {p\over \pi (s-1)} 
\int_0^{\pi/2}  
\big\{D_p(\theta)^{s-1} - D_p(\pi-\theta)^{s-1}\big\} 
\cos\theta \d \theta.
\cr}$$
This implies that
$$\bigg|{E_p(s)\over E_p(\sigma)}\bigg|
= \bigg|{\sigma-1\over s-1}\bigg|
\bigg|{E_p^*(s)\over E_p^*(\sigma)}\bigg|
\leqno(5.2)$$
with
$$E_p^*(s)
:= \int_0^{\pi/2} 
\big\{D_p(\theta)^{s-1} - D_p(\pi-\theta)^{s-1}\big\} 
\cos\theta \d \theta.$$

\medskip

$1^\circ$
Case of $\sigma^{1/\delta}<|\tau|\le y^{1/\delta}$

\smallskip

Write
$$E_p^*(s)
= \int_0^{\pi/2} D_p(\theta)^{s-1}
\big\{1 - \Delta_p(\theta)^{s-1}\big\} 
\cos\theta \d \theta$$
with
$$\Delta_p(\theta)
:= {1-2p^{-1}\cos\theta+p^{-2}\over 1+2p^{-1}\cos\theta+p^{-2}}.$$
It is clear that for all $p$,
the function $\theta\mapsto\Delta_p(\theta)$ 
is increasing on $[0, \pi/2]$.
It follows that
$$\eqalign{E_p^*(\sigma)
& \ge \int_0^{\pi/4} D_p(\theta)^{\sigma-1} 
\big\{1 - \Delta_p(\theta)^{\sigma-1}\big\} 
\cos\theta \d \theta
\cr
& \ge \big\{1 - \Delta_p(\pi/4)^{\sigma-1}\big\}
\int_0^{\pi/4} D_p(\theta)^{\sigma-1}  
\cos\theta \d \theta
\cr}$$
for all $p$ and $\sigma\ge 1$.
This implies that
$$\bigg|{1\over E_p^*(\sigma)}
\int_0^{\pi/4} D_p(\theta)^{\sigma-1} 
\cos\theta \d \theta\bigg|
\le {1\over 1 - \Delta_p(\pi/4)^{\sigma-1}}.
\leqno(5.3)$$
Similarly since the function 
$\theta\mapsto D_p(\theta)^{\sigma-1} \cos\theta$
is decreasing on $[0, \pi/2]$ for all $p$ and $\sigma\ge 2$,
we can deduce, via (5.3), that
$$\bigg|{1\over E_p^*(\sigma)}
\int_{\pi/4}^{\pi/2} D_p(\theta)^{\sigma-1} 
\cos\theta \d \theta\bigg|
\le {1\over 1 - \Delta_p(\pi/4)^{\sigma-1}}.
\leqno(5.4)$$
From (5.3) and (5.4), we deduce that
$$\bigg|{E_p^*(s)\over E_p^*(\sigma)}\bigg|
\le {2\over 1 - \Delta_p(\pi/4)^{\sigma-1}}.$$
It is easy to verify that
for all $p\ge \sigma\ge 2$, we have
$$\Delta_p\bigg({\pi\over 4}\bigg)^{\sigma-1}
\le \bigg(1-{\sqrt 2\over p}+{1\over p^2}\bigg)^{\sigma-1}
\le 1 - {\sigma-1\over 4p}.$$
Combining these estimates with (5.2), we obtain
$$\bigg|{E_p(s)\over E_p(\sigma)}\bigg|
\le {8p\over |s-1|}
\le {p^4\over |\tau|}
\qquad(p\ge \sigma).$$
By multiplying this inequality for $\sigma<p\le |\tau|^\delta\;(\le y)$
and the trivial inequality $|E_p(s)|\le |E_p(\sigma)|$ 
for the others $p$,
we deduce, via the prime number theorem, that
$$\eqalign{\bigg|{E(s, y)\over E(\sigma, y)}\bigg|
& \le \exp\bigg\{-\sum_{\sigma<p\le |\tau|^\delta} \log|\tau| 
+ 4 \sum_{\sigma<p\le |\tau|^\delta} \log p\bigg\}
\cr
& \le e^{-\{1/\delta-4+o(1)\}|\tau|^\delta}.
\cr}$$

\smallskip

$2^\circ$
Case of $c_1\sigma^{1/2}\log\sigma\le |\tau|\le \sigma^{1/\delta}$

\smallskip

For $p\ge \sigma^{1/2}\ge 2$, we can write
$$\eqalign{|E_p^*(s)|
& \le \int_0^{\pi/2} 
\big\{D_p(\theta)^{\sigma-1} 
+ D_p(\pi-\theta)^{\sigma-1}\big\} 
\cos\theta \d \theta
\cr
& = \bigg\{1+O\bigg({\sigma\over p^2}\bigg)\bigg\}
\int_0^{\pi/2} 
\big(e^{2[(\sigma-1)/p]\cos\theta}
+e^{-2[(\sigma-1)/p]\cos\theta}\big) 
\cos\theta \d \theta
\cr}$$
and
$$\eqalign{|E_p^*(\sigma)|
& = \int_0^{\pi/2} 
\big\{D_p(\theta)^{\sigma-1} 
- D_p(\pi-\theta)^{\sigma-1}\big\} 
\cos\theta \d \theta
\cr
& = \bigg\{1+O\bigg({\sigma\over p^2}\bigg)\bigg\}
\int_0^{\pi/2} 
\big(e^{2[(\sigma-1)/p]\cos\theta}
-e^{-2[(\sigma-1)/p]\cos\theta}\big) 
\cos\theta \d \theta.
\cr}$$ 
From these, we deduce that
$$\bigg|{E_p^*(s)\over E_p^*(\sigma)}\bigg|
\le \bigg\{1+O\bigg({\sigma\over p^2}+{1\over e^{\sigma/p}}\bigg)\bigg\}
\qquad(2\le \sigma^{1/2}\le p\le \sigma)
\leqno(5.5)$$
where we have used the following facts
$$\int_0^{\pi/2} e^{2[(\sigma-1)/p]\cos\theta} \cos\theta \d \theta
\gg e^{\sigma/p}
\qquad{\rm and}\qquad
\int_0^{\pi/2} e^{-2[(\sigma-1)/p]\cos\theta} \cos\theta \d \theta
\ll 1.$$
Inserting (5.5) into (5.2), 
for $2\le \sigma^{1/2}\le p\le \sigma$ we obtain 
$$\eqalign{\bigg|{E_p(s)\over E_p(\sigma)}\bigg|
& \le \exp\bigg\{-\log\bigg|{s-1\over \sigma-1}\bigg|
+ C\bigg({\sigma\over p^2}+{1\over e^{\sigma/p}}\bigg)\bigg\}
\cr\noalign{\medskip}
& \le \cases{
e^{-\tau^2/(2\sigma^2)+C\sigma/p^2+Ce^{-\sigma/p}}
& if $\;\;3\le |\tau|\le \sigma$,
\cr\noalign{\medskip}
e^{-{1\over 2}\log(1+\tau^2/\sigma^2) + C\sigma/p^2+Ce^{-\sigma/p}}
& if $\;\;\sigma\le |\tau|\le \sigma^{1/\delta}$,
\cr}
\cr}$$
where $C>0$ is an absolute constant.

Now by multiplying these inequalities 
for $\sigma/(4\log\sigma)\le p\le \sigma/(2\log\sigma)$ and 
the trivial inequality 
$|E_p(s)|\le E_p(\sigma)$ for the other $p$,
we get
$$\leqalignno{\bigg|{E(s, y)\over E(\sigma, y)}\bigg|
& \le \exp\bigg\{-
\sum_{\sigma/(4\log\sigma)\le p\le \sigma/(2\log\sigma)} 
\bigg({\tau^2\over 2\sigma^2}
- {C\sigma\over p^2}-{C\over e^{\sigma/p}}\bigg)\bigg\}
\cr
& \le \exp\bigg\{-
\bigg({\tau^2\over 16\sigma(\log\sigma)^2}
- 10C
-{10C\over \sigma\log\sigma}\bigg)
\bigg\}
\cr
& \le \exp\bigg\{-{c_2\tau^2\over \sigma(\log\sigma)^2}\bigg\}
\cr}$$
if $c_1\sigma^{1/2}\log\sigma\le |\tau|\le \sigma$,
and
$$\leqalignno{\bigg|{E(s, y)\over E(\sigma, y)}\bigg|
& \le \exp\bigg\{-
\sum_{\sigma/(4\log\sigma)\le p\le \sigma/(2\log\sigma)} 
\bigg[{1\over 2}\log\bigg(1+{\tau^2\over \sigma^2}\bigg)
- {C\sigma\over p^2}-{C\over e^{\sigma/p}}\bigg]\bigg\}
& (5.6)
\cr
& \le \exp\bigg\{-
\bigg[
{\sigma\over 8\log\sigma}\log\bigg(1+{\tau^2\over \sigma^2}\bigg)
- 10C
-{10C\over \sigma\log\sigma}\bigg]
\bigg\}
\cr\noalign{\medskip}
& \le \exp\big\{-c_3|\tau|^\delta\big\}
\cr}$$
if $\sigma\le |\tau|\le \sigma^{1/\delta}$.
This completes the proof.
\hfill
$\square$

\goodbreak
\vskip 8mm

\centerline{\pptitre \S\ 6. Proof of Theorem 3}

\bigskip

We follow the argument of Granville \& Soundararajan [4]
to prove Theorem 3.
We shall divide the proof in several steps which are embodied 
in the following lemmas. 

The first one is a classic integration formula (see [4], page 1019).

\proclaim Lemma 6.1.
Let $c>0$, $\lambda>0$ and $N\in \N$.
Then we have
$${1\over 2\pi i}\int_{c-i\infty}^{c+i\infty}
y^s \bigg({e^{\lambda s}-1\over \lambda s}\bigg)^N {\d s\over s}
= \cases{
0         & if $\;0<y<e^{-\lambda N}$,
\cr\noalign{\vskip 2mm}
\in [0, 1] & if $\;e^{-\lambda N}\le y<1$,
\cr\noalign{\vskip 2mm}
1         & if $\;y\ge 1$.
\cr}
\leqno(6.1)$$

\smallskip

The second one is an analogue for (3.6) and (3.7) of [4]
(see also Lemma 3.1 of [20]).

\proclaim Lemma 6.2.
Let $t\ge 1$, $y\ge 2e^t$ and $0<\lambda\le e^{-t}$.
Then we have
$$\displaylines{
\rlap{\rm (6.2)}
\hfill
\Phi(t, y) 
\le {1\over 2\pi i}\int_{\kappa-i\infty}^{\kappa+i\infty}
{E(s, y)\over (e^\gamma t)^{2s}}
{e^{\lambda s}-1\over \lambda s} {\d s \over s} 
\le \Phi(te^{-\lambda}, y),
\hfill
\cr\noalign{\smallskip}
\rlap{\rm (6.3)}
\hfill
\Phi(te^{-\lambda}, y) - \Phi(t, y)   
\le {1\over 2\pi i}\int_{\kappa-i\infty}^{\kappa+i\infty}
{E(s, y)\over (e^\gamma t)^{2s}} 
{e^{\lambda s}-1\over \lambda s} 
\big(e^{2\lambda s} - e^{-2\lambda s}\big) {\d s \over s}.
\hfill
\cr}$$

\noindent{\sl Proof}.
Denote by ${\bf 1}_X(\omega)$
the characteristic function of the set $X\subset \Omega$.
Then by Lemma 6.1 with $N=1$ and $c=\kappa$, we have
$${\bf 1}_{\{\omega\in \Omega : 
L(1, g^\natural(\omega); y)>(e^\gamma t)^2\}}(\omega)
\le {1\over 2\pi i}\int_{\kappa-i\infty}^{\kappa+i\infty}
\bigg({L(1, g^\natural(\omega); y)\over (e^\gamma t)^2}\bigg)^s
{e^{\lambda s}-1\over \lambda s^2} \d s.$$
Integrating over $\Omega$ and 
interchanging the order of integrations yield
$$\eqalign{
\Phi(t, y) 
& \le \int_{\Omega}
\bigg\{{1\over 2\pi i}\int_{\kappa-i\infty}^{\kappa+i\infty}
\bigg({L(1, g^\natural(\omega); y)\over (e^\gamma t)^2}\bigg)^s
{e^{\lambda s}-1\over \lambda s^2} \d s\bigg\}
\d \mu(\omega)
\cr\noalign{\smallskip}
& = {1\over 2\pi i}\int_{\kappa-i\infty}^{\kappa+i\infty}
{E(s, y)\over (e^\gamma t)^{2s}}
{e^{\lambda s}-1\over \lambda s^2} \d s. 
\cr}$$
This proves the first inequality of (6.2).
The second can be treated by noticing that
$$\eqalign{
{\bf 1}_{\{\omega\in \Omega : 
L(1, g^\natural(\omega); y)>(e^{\gamma-\lambda}t)^2\}}(\omega) 
& = {\bf 1}_{\{\omega\in \Omega : 
L(1, g^\natural(\omega); y)>(e^{\gamma}t)^2\}}(\omega) 
\cr\noalign{\smallskip}
& \, + {\bf 1}_{\{\omega\in \Omega : 
(e^{\gamma}t)^2
\ge L(1, g^\natural(\omega); y)
>(e^{\gamma-\lambda}t)^2\}}(\omega) 
\cr\noalign{\smallskip}
& \ge {1\over 2\pi i}\int_{\kappa-i\infty}^{\kappa+i\infty}
\bigg({L(1, g^\natural(\omega); y)\over (e^\gamma t)^2}\bigg)^s
{e^{\lambda s}-1\over \lambda s^2} \d s.
\cr}$$
From (6.2), we can deduce
$$\eqalign{
\Phi(te^{-\lambda}, y) - \Phi(t, y)   
& \le {1\over 2\pi i}\int_{\kappa-i\infty}^{\kappa+i\infty}
{E(s, y)\over (e^{\gamma-\lambda}t)^{2s}} 
{e^{\lambda s}-1\over \lambda s^2} \d s 
- {1\over 2\pi i}\int_{\kappa-i\infty}^{\kappa+i\infty}
{E(s, y)\over (e^{\gamma+\lambda} t)^{2s}} 
{e^{\lambda s}-1\over \lambda s^2} \d s
\cr
& = {1\over 2\pi i}\int_{\kappa-i\infty}^{\kappa+i\infty}
{E(s, y)\over (e^{\gamma} t)^{2s}} 
{e^{\lambda s}-1\over \lambda s^2} 
\big(e^{2\lambda s} - e^{-2\lambda s}\big) \d s.
\cr}$$
This completes the proof.
\hfill
$\square$

\proclaim Lemma 6.3.
Let $t\ge 1$, $y\ge 2e^t$ and $0<\kappa\lambda\le 1$.
Then we have
$${1\over 2\pi i}\int_{\kappa-i\kappa}^{\kappa+i\kappa} 
{E(s, y)\over (e^{\gamma} t)^{2s}} 
{e^{\lambda s}-1\over \lambda s^2} \d s
= {E(\kappa, y)\over \kappa\sqrt{2\pi\sigma_2}(e^\gamma t)^{2\kappa}}
\bigg\{1 +
O\bigg(\kappa \lambda + {\log\kappa\over \kappa}\bigg)\bigg\}.$$

\noindent{\sl Proof}.
First in view of (4.24) we write, for $s=\kappa+i\tau$ 
and $|\tau|\le \kappa$,
$$E(s, y)
= \exp\bigg\{\sigma_0 
+ i\sigma_1\tau  
- {\sigma_2\over 2}\tau^2 
- i{\sigma_3\over 6}\tau^3 
+ O\big(\sigma_4\tau^4\big)\bigg\}$$
and
$${e^{\lambda s}-1\over \lambda s^2}
= {1\over \kappa} 
\bigg\{1 - {i\over \kappa}\tau 
+ O\bigg(\kappa\lambda + {\tau^2\over \kappa^2}\bigg)\bigg\}.$$
Since
$\sigma_1=\log t+\gamma$,
we have
$$\eqalign{{E(s, y)\over (e^{\gamma} t)^{2s}} 
{e^{\lambda s}-1\over \lambda s^2}
& = {E(\kappa, y)\over \kappa (e^\gamma t)^{2\kappa}}
e^{-(\sigma_2/2) \tau^2}
\bigg\{1 - {i\over \kappa}\tau - i{\sigma_3\over 6}\tau^3
+ O\big(R(\tau)\big)
\bigg\}
\cr}$$
with
$$R(\tau) := \kappa\lambda 
+ \kappa^{-2} \tau^2 
+ \sigma_4 \tau^4
+ \sigma_3^2 \tau^6.$$
Now we integrate the last expression over $|\tau|\le \kappa$ 
to obtain
$${1\over 2\pi i}\int_{\kappa-i\kappa}^{\kappa+i\kappa} 
{E(s, y)\over (e^{\gamma} t)^{2s}} 
{e^{\lambda s}-1\over \lambda s^2} \d s
= {E(\kappa, y)\over 2\pi \kappa (e^\gamma t)^{2\kappa}}
\int_{-\kappa}^{\kappa} e^{-(\sigma_2/2) \tau^2} 
\big\{1 + O\big(R(\tau)\big)\big\} \d\tau,
\leqno(6.4)$$
where we have used the fact that
the integrals involving $(i/\kappa)\tau$ 
and $(i\sigma_3/6)\tau^3$ vanish.

On the other hand, using 
lemmas 4.5 and 4.6
we have
$$\eqalign{
\int_{-\kappa}^{\kappa} 
e^{-(\sigma_2/2) \tau^2} \d\tau
& = \sqrt{2\pi\over \sigma_2}
\bigg\{1 
+ O\bigg(\exp\bigg\{-{1\over 2}\kappa^2\sigma_2\bigg\}\bigg)\bigg\},
\cr\noalign{\smallskip}
\int_{\kappa-i\kappa}^{\kappa+i\kappa} 
e^{-(\sigma_2/2) \tau^2} R(\tau) \d\tau
& \ll {1\over \sqrt{\sigma_2}}
\bigg(\kappa\lambda 
+ {1\over \kappa^2 \sigma_2} 
+ {\sigma_3^2\over \sigma_2^3} 
+ {\sigma_4\over \sigma_2^2}
\bigg)
\cr
& \ll {1\over \sqrt{\sigma_2}}
\bigg(\kappa\lambda + {\log\kappa\over \kappa}\bigg).
\cr}$$
Inserting these into (6.4), we obtain the desired result.
\hfill
$\square$

\proclaim Lemma 6.4.
Let $\delta$ and $c_3$ be two constants determined by Lemma 5.1.
Then we have
$$\leqalignno{\qquad
\int_{\kappa\pm i\kappa}^{\kappa\pm i\infty}
{E(s, y)\over (e^{\gamma} t)^{2s}} 
{e^{\lambda s}-1\over \lambda s^2} \d s 
& \ll {E(\kappa, y)\over \kappa\sqrt{\sigma_2} (e^\gamma t)^{2\kappa}}
R_1,
& (6.5)
\cr\noalign{\medskip}
\int_{\kappa-i\infty}^{\kappa+i\infty}
{E(s, y)\over (e^{\gamma} t)^{2s}} 
{e^{\lambda s}-1\over \lambda s^2} 
\big(e^{2\lambda s} - e^{-2\lambda s}\big) \d s
& \ll {E(\kappa, y)\over \kappa\sqrt{\sigma_2}
(e^\gamma t)^{2\kappa}} R_2,
& (6.6)
\cr}$$
uniformly for $t\ge 1$, $y\ge 2e^t$, $\kappa\ge 2$ 
and $0<\lambda\kappa\le 1$, 
where
$$\eqalign{
R_1
& := \lambda^{-1} e^{-c_3\kappa^\delta} 
+ \lambda^{-1} (\kappa/\log\kappa)^{1/2} y^{-1/\delta},
\cr\noalign{\smallskip}
R_2
& := \lambda \kappa(\log\kappa)^{1/2} 
+ e^{-(c_3/2)\kappa^\delta}
+ \lambda^{-1} (\kappa/\log\kappa)^{1/2} y^{-1/\delta}.
\cr}$$

\noindent{\sl Proof}.
We split the integral in (6.5) into two parts according to
$\kappa\le |\tau|\le y^{1/\delta}$
or
$|\tau|\ge y^{1/\delta}$.
Using Lemma 5.1 with $\sigma=\kappa$ and 
the inequality $(e^{\lambda s}-1)/s^2\ll 1/\tau^2$, 
the integral in (6.5) is
$$\ll {E(\kappa, y)\over (e^\gamma t)^{2\kappa} \lambda}
\bigg({e^{-c_3\kappa^\delta}\over \kappa}
+ {1\over y^{1/\delta}}
\bigg),$$
which implies (6.5), in view of Lemma 4.5 with $J=1$.

\smallskip

Similarly we split the integral in (6.6) 
into four parts according to 
$$|\tau|\le c_1\kappa^{1/2}\log\kappa,
\qquad 
c_1\kappa^{1/2}\log\kappa<|\tau|\le \kappa,
\qquad 
\kappa<|\tau|\le y^{1/\delta},
\qquad
|\tau|\ge y^{1/\delta}.$$
By Lemma 5.1 with $\sigma=\kappa$ and the inequalities
$$\eqalign{(e^{\lambda s}-1)/\lambda s
& \ll \min\{1, \; 1/(\lambda |\tau|)\},
\cr\noalign{\vskip 1mm}
(e^{2\lambda s} - e^{-2\lambda s})/s
& \ll \min\{\lambda, \; 1/|\tau|\},
\cr}$$
the integral in (6.6) is, as before, 
$$\ll_\varepsilon {E(\kappa, y)\over (e^\gamma t)^{2\kappa}}
\Big(\lambda\kappa^{1/2}\log\kappa
+ e^{-c_3\kappa^\delta}
+ \lambda^{-1} y^{-1/\delta}
\Big),$$
which implies (6.6), as before.
\hfill
$\square$

\medskip

Now we are ready to complete the proof of Theorem 3.
Lemma 6.3 and (6.5) of Lemma 6.4 give 
$${1\over 2\pi i}\int_{\kappa-i\infty}^{\kappa+i\infty}
{E(s, y)\over (e^{\gamma} t)^{2s}} 
{e^{\lambda s}-1\over \lambda s^2} \d s  
= {E(\kappa, y)\over \kappa\sqrt{2\pi\sigma_2}
(e^\gamma t)^{2\kappa}}
\big\{1 + O\big(R'\big)\big\}
\leqno(6.7)$$
where
$$R' 
:= {\log\kappa\over \kappa}
+ \kappa\lambda 
+ {e^{-c_3\kappa^\delta}+(\kappa/\log\kappa)^{1/2} y^{-1/\delta}
\over \lambda}.
$$
Taking $\lambda = \kappa^{-2}$ 
and noticing $y\ge 2e^t\asymp \kappa$ and $1/\delta>4$, 
we deduce
$$R'\ll {t/e^t}.
\leqno(6.8)$$
Combining (6.7) and (6.8) with (6.2), we obtain
$$\Phi(t, y)
\le {E(\kappa, y)\over \kappa\sqrt{2\pi\sigma_2}
(e^\gamma t)^{2\kappa}}
\bigg\{1 + O\bigg({t\over e^t}\bigg)\bigg\}
\le \Phi(te^{-\lambda}, y)
\leqno(6.9)$$   
uniformly for $t\ge 1$, $y\ge 2e^t$ and $0<\lambda\le e^{-t}$.

On the other hand, 
(6.3) of Lemma 6.2 and (6.6) of Lemma 6.4 imply
$$\eqalign{\Phi(te^{-\lambda}, y) - \Phi(t, y)   
& \ll {E(\kappa, y)\over \kappa\sqrt{\sigma_2} (e^\gamma t)^{2\kappa}}
\bigg(\lambda \kappa(\log\kappa)^{1/2}
+ {(\kappa/\log\kappa)^{1/2}\over e^{c_3\kappa^\delta}}
+ {(\kappa/\log\kappa)^{1/2}\over \lambda y^{1/\delta}}\bigg)
\cr
& \ll {E(\kappa, y)\over \kappa\sqrt{\sigma_2}
(e^\gamma t)^{2\kappa}}
\bigg(\lambda \kappa(\log\kappa)^{1/2}
+ {(\kappa/\log\kappa)^{1/2}\over e^{c_3\kappa^\delta}}\bigg)
\cr}$$
when $y^{-1/(2\delta)}\kappa^{-1/2}(\log\kappa)^{-1}
\le \lambda\le \kappa^{-1}$.
Since $\Phi(te^{-\lambda}, y) - \Phi(t, y)$ 
is a non-decreasing function of $\lambda$,
we deduce 
$$\Phi(te^{-\lambda}, y) - \Phi(t, y)   
\ll {E(\kappa, y)\over \kappa\sqrt{\sigma_2} (e^\gamma t)^{2\kappa}}
\bigg(\lambda \kappa(\log\kappa)^{1/2}
+ {(\kappa/\log\kappa)^{1/2}\over e^{c_3\kappa^\delta}}
+ {\kappa(\log\kappa)^{1/2}\over y^{1/(2\delta)}}\bigg)
\leqno(6.10)$$
uniformly for $t\ge 1$, $y\ge 2e^t$ and $0<\lambda\le e^{-t}$.
Obviously the estimates (6.9) and (6.10) imply 
the desired result.
This completes the proof of Theorem 3. 
\hfill
$\square$

\vskip 8mm

\centerline{\pptitre \S\ 7. Proof of Theorem 4}

\bigskip

Using Lemmas 4.1 and 4.5, we can write
$$\eqalign{
{E(\kappa, y)\over \kappa\sqrt{2\pi\sigma_2}(e^\gamma t)^{2\kappa}}
& = \exp\bigg\{\phi(\kappa, y)
- 2\kappa (\gamma + \log t)
+ O(\log\kappa)\bigg\}
\cr
& = \exp\bigg\{\kappa \bigg(
2\log_2\kappa
- 2\log t
+ \sum_{j=1}^J {b_{j, 0}\over (\log\kappa)^j}
+ O_J\big(R_J(\kappa, y)\big)
\bigg)
\bigg\}.
\cr}$$
On the other hand, 
Lemma 4.3 and (1.17) imply that
$$2\log_2\kappa + 2\gamma 
+ \sum_{j=1}^J {b_{j, 1}\over (\log\kappa)^j}
+ O_J\big(R_J(\kappa, y)\big)
= 2(\log t + \gamma).$$
Combining these estimates, we can obtain
$$\eqalign{
{E(\kappa, y)\over \kappa\sqrt{2\pi\sigma_2}(e^\gamma t)^{2\kappa}}
& = \exp\bigg\{-\kappa
\bigg[\sum_{j=1}^J {b_{j, 1}-b_{j, 0}\over (\log\kappa)^j}
+ O_J\big(R_J(\kappa, y)\big)\bigg]
\bigg\}.
\cr}$$
In view of (1.21), (4.2) and (4.15),
we have $b_{j, 1} - b_{j, 0} = a_j$.
This completes the proof.
\hfill
$\square$

\vskip 8mm

\centerline{\pptitre \S\ 8. Proof of Corollary 5}

\bigskip

We first prove an asymptotic developpment of $\kappa(t, y)$
in $t$.

\proclaim Lemma 8.1.
For each integer $J\ge 1$, 
there are computable constants $\gamma_0, \gamma_1, \dots, \gamma_J$ 
such that the asymptotic formula
$$\kappa(t, y)
= e^{t-\gamma_0}
\bigg\{1 + \sum_{j=1}^J {\gamma_j\over t^j} 
+ O_J\big(R_J^*(t, y)\big)
\bigg\}
\leqno(8.1)$$
holds uniformly for $t\ge 1$ and $y\ge 2e^t$,
where 
$$R_N^*(t, y)
:= {1\over t^{N+1}}+{e^tt\over y\log y}.$$
Further $\gamma_0$ is given by (1.24) and  
$\gamma_1 = - {1\over 8} b_{1,1}^2 - {1\over 4} b_{2,1}$.

\noindent{\sl Proof}.
By Lemma 4.3 and (1.17), we have
$$2\log t 
= 2\log_2\kappa + 
\sum_{j=1}^{J+1} {b_{j, 1}\over (\log\kappa)^j}
+ O_J\big(R_{J+1}(\kappa, y)\big),
\leqno(8.2)$$
where $R_J(\kappa, y)$ is defined as in (1.20). 
From (8.2), we easily deduce that
$$\eqalign{t 
& = (\log\kappa) \prod_{j=1}^{J+1}
\exp\bigg\{{b_{j,1}\over 2(\log\kappa)^j}\bigg\}
\exp\big\{O_J\big(R_{J+1}(\kappa, y)\big)\big\}
\cr
& = (\log\kappa) \prod_{j=1}^{J+1}
\bigg\{\sum_{m_j=0}^{J+1} {1\over m_j!}
\bigg({b_{j,1}\over 2(\log\kappa)^j}\bigg)^{m_j}
+ O_J\big(R_{J+1}(\kappa, y)\big)\bigg\}.
\cr}$$
Developping the product, we get
$$\eqalign{t 
& = (\log\kappa) 
\bigg\{\sum_{j=0}^{J+1} {b_j'\over (\log\kappa)^j}
+ O_J\big(R_{J+1}(\kappa, y)\big)\bigg\},
\cr}$$
where
$$
\eqalign{
b_j' 
& :=
\sum_{\scriptstyle m_1\ge 0, \dots, m_{J+1}\ge 0
\atop\scriptstyle m_1+2m_2+\cdots+(J+1)m_{J+1}=j}
{b_{1,1}^{m_1}\cdots b_{J+1,1}^{m_{J+1}}\over 
(2m_1)!!\cdots (2m_{J+1})!!}
\cr
& \,=
\sum_{\scriptstyle m_1\ge 0, \dots, m_{j}\ge 0
\atop\scriptstyle m_1+2m_2+\cdots+jm_{j}=j}
{b_{1,1}^{m_1}\cdots b_{j,1}^{m_{j}}\over
(2m_1)!!\cdots (2m_{j})!!}.
\cr}$$
Since $b_0'=1$ and $b_1'=b_{1,1}/2=\gamma_0$,
the preceeding asymptotic formula can be written as
$$t 
= \log\kappa + \gamma_0
+ \sum_{j=1}^J {b_{j+1}'\over (\log\kappa)^j}
+ O_J\big(R_J^*(t, y)\big),
\leqno(8.3)$$
where we have used the fact that $\kappa(t, y)\asymp e^t$
(see Lemma 2.3)
and $(\log k)R_{J+1}(\kappa, y)\asymp R_J^*(t, y)$.

With the help of (8.3), 
a simple recurrence argument shows that 
there are constants $\gamma'_n$ such that
$$t 
= \log\kappa + \sum_{j=0}^J {\gamma'_j\over t^j}
+ O_J\big(R_J^*(t, y)\big).
\leqno(8.4)$$
In fact taking $J=0$ in (8.3), 
we see that (8.4) holds for $J=0$.
Suppose that it holds for $0, \dots, J-1$, i.e.
$$t 
= \log\kappa + \sum_{i=0}^{J-j-1} {\gamma'_i\over t^i}
+ O\big(R_{J-j-1}^*(t, y)\big)
\qquad(j=0, \dots, J-1),$$
which is equivalent to
$$\log\kappa 
= t \bigg\{1 - \sum_{i=1}^{J-j} {\gamma'_{i-1}\over t^i}
+ O\bigg({R_{J-j-1}^*(t, y)\over t}\bigg)\bigg\}
\qquad(j=0, \dots, J-1).
\leqno(8.5)$$
This holds also for $j=J$ 
if we use the convention: 
$$\sum_{i=0}^{-1}=0
\qquad{\rm and}\qquad
R_{-1}^*(t, y) := 1,$$
since $\log\kappa=t + O(1)$.
Inserting it into (8.3), 
we easily see that (8.4) holds also for $J$.
In particular we have 
$$\textstyle
\gamma'_1 = b_2' = {1\over 8}b_{1,1}^2 + {1\over 4}b_{2,1}.$$

Now (8.1) is an immediate consequence of (8.4) with
$$\gamma_j :=
\sum_{\scriptstyle m_1\ge 0, \dots, m_J\ge 0
\atop\scriptstyle m_1+2m_2+\cdots+Jm_J=j}
(-1)^{m_1+\cdots+m_J}
{{\gamma_1'}^{m_1}\cdots {\gamma_J'}^{m_J}\over m_1!\cdots m_J!}.$$
This completes the proof.
\hfill
$\square$

\medskip

Now we are ready to prove Corollary 5.

Using (8.5), we have
$$\leqalignno{\sum_{j=1}^J {a_j\over (\log\kappa)^j}
& = \sum_{j=1}^J {a_j\over t^j}
\bigg\{1 - \sum_{i=1}^{J-j} {\gamma'_{i-1}\over t^i} 
+ O_N\bigg({R_{J-j-1}^*(t, y)\over t}\bigg)\bigg\}^{-j}
& (8.6)
\cr
& = \sum_{j=1}^J {\rho_j\over t^j}
+ O_J\bigg({R_{J-2}^*(t, y)\over t^2}\bigg),
\cr}$$
where the $\rho_n$ are constants. 
In particular we have
$\rho_1 = a_1 = 1$
and
$\rho_2 = \gamma_0 + a_2$.

Now Theorem 4, (8.1) and (8.6) imply the result of Corollary with
$$a^*_1
= \rho_1 = 1,
\qquad
a^*_j
= \rho_j + \sum_{i=1}^{j-1} \gamma_i\rho_{j-i}
\quad(j\ge 2).$$
This completes the proof of Corollary 5.
\hfill
$\square$

\vskip 10mm

\centerline{\pptitre \S\ 9. Proof of Theorem 2}

\bigskip

For each $\eta\in (0, \dm)$, define
$${\rm H}_k^+(1; \eta)
:= \big\{f\in {\rm H}_k^*(1) : 
L(s, f)\not=0, \, s\in {\cal S}\big\},$$
where 
${\cal S} 
:= \{s:=\sigma+i\tau :\, \sigma\ge 1-\eta, \, |\tau|\le 100k^\eta\}\cup
\{s:=\sigma+i\tau:\, \sigma\ge 1,\, \tau\in\R\}$, and
$${\rm H}_k^-(1; \eta)
:= {\rm H}_k^*(1)\sset {\rm H}_k^+(1; \eta).
$$
Then we have (see [10], (1.11))
$$\big|{\rm H}_k^-(1; \eta)\big|
\ll_\eta k^{31\eta}.
\leqno(9.1)$$

Our starting point in the proof of Theorem 2 is the
evaluation of the moments of $L(1,f)$. For this, we recall a 
particular case of Proposition 6.1 of [10].

\proclaim Lemma 9.1.
Let $\eta\in (0, {1\over 31})$ be fixed.
There are two positive constants  
$c_i=c_i(\eta)\;(i=4, 5)$ such that
$$\sum_{f\in {\rm H}_k^{+}(1; \eta)} \omega_f L(1, f)^s
= E(s) + O_\eta\big(e^{-c_4\log k/\log_2k}\big)
\leqno(9.2)$$
uniformly for
$$k\ge 16,
\qquad
2\mid k
\qquad\hbox{and}\qquad
|s|\le 2T_k
\leqno(9.3)$$
with $$T_k:=c_5\log k/(\log_2k\log_3k).$$
Here $E(s)$ is defined by (1.15).

Let $\kappa(t, y)$ be the saddle-point determined by (1.17)
and $\kappa_t:=\kappa(t, \infty)$.
For $k\ge 16, 2\mid k$, $\lambda>0$, $N\in \N$ and $t>0$,
introduce the two integrals
$$I_1(k, t; \lambda, N)
:={1\over 2\pi i}\int_{\kappa_t-i\infty}^{\kappa_t+i\infty}
\sum_{f\in {\rm H}_k^+(1;\eta)}\omega_f
\bigg({L(1,f)\over (e^\gamma t)^{2}}\bigg)^s 
\bigg({e^{\lambda s}-1\over \lambda s}\bigg)^{2N} {\d s\over s}$$
and
$$I_2(k, t; \lambda, N)
:={1\over 2\pi i}\int_{\kappa_t-i\infty}^{\kappa_t+i\infty} 
{E(s)\over (e^\gamma t)^{2s}} 
\bigg({e^{\lambda s}-1\over \lambda s}\bigg)^{2N} {\d s\over s}.$$

\proclaim Lemma 9.2.
Let $\eta\in (0, {1\over 200}]$ be fixed.
Then we have
$$\leqalignno{
\widetilde{F}_k(t)+O_\eta\big(k^{-5/6}\big)
& \le I_1(k, t; \lambda, N)
\le \widetilde{F}_k(te^{-\lambda N})+O_\eta\big(k^{-5/6}\big),
& (9.4)
\cr\noalign{\vskip 1mm}
\Phi(t)
& \le I_2(k, t; \lambda, N)
\le \Phi(te^{-\lambda N})
& (9.5)
\cr}$$
uniformly for $k\ge 16$, $2\mid k$, $\lambda>0$, $N\in \N$ and $t>0$.
The implied constants depend on $\eta$ only.

\noindent{\sl Proof}.
By exchanging the order of sommation
and by using Lemma 6.1 with $c=\kappa_t$, we obtain
$$\eqalign{
I_1(k, t; \lambda, N)
& = \sum_{f\in {\rm H}_k^+(1;\eta)}
{\omega_f\over 2\pi i}\int_{\kappa_t-i\infty}^{\kappa_t+i\infty}
\bigg({L(1,f)\over (e^\gamma t)^{2}}\bigg)^s 
\bigg({e^{\lambda s}-1\over \lambda s}\bigg)^{2N} {\d s\over s},
\cr 
& \ge \sum_{f\in {\rm H}_k^+(1;\eta), \;
L(1,f)\ge (e^\gamma t)^{2}} \omega_f.
\cr}$$
In view of the second estimate of (1.7) and of (9.1), 
we reintroduce the missing forms
$$\eqalign{I_1(k, t; \lambda, N)
& \ge \sum_{f\in {\rm H}_k^*(1),\; L(1, f)\ge (e^\gamma t)^2} 
\omega_f
+ O\Big(\sum_{f\in {\rm H}_k^*\sset{\rm H}_k^+(1;\eta)}\omega_f\Big)
\cr
& \ge \sum_{f\in {\rm H}_k^*(1),\; L(1, f)\ge (e^\gamma t)^2} 
\omega_f
+ O\big(k^{-1+31\eta}\log k\big).
\cr}$$
Clearly this implies the first inequality of (9.4),
thanks to (1.6) and (1.7).

Similarly, using Lemma 6.1 with $c=\kappa_t$, we find
$$\eqalign{I_1(k, t; \lambda, N)
& \le \sum_{\scriptstyle f\in {\rm H}_k^+(1;\eta)
\atop\scriptstyle L(1,f)\ge (e^\gamma t)^{2}} \omega_f
+ \sum_{\scriptstyle f\in {\rm H}_k^+(1;\eta)
\atop\scriptstyle 
(e^\gamma te^{-\lambda N})^{2}
\le L(1,f)<(e^\gamma t)^{2}} \omega_f
\cr
& = \sum_{\scriptstyle f\in {\rm H}_k^+(1;\eta)
\atop\scriptstyle 
L(1,f)\ge (e^\gamma te^{-\lambda N})^{2}} \omega_f.
\cr}$$
As before, we can easily show that 
the last sum is 
$\le \widetilde{F}_k(te^{-\lambda N})+O\big(k^{-5/6}\big)$.

The estimates (9.5) can be proved in the same way as (6.2).
\hfill
$\square$

\proclaim Lemma 9.3.
Let $\eta\in (0, {1\over 200}]$ be fixed and
$c_4$ be the positive constant given by Lemma 9.1.
Then we have
$$\eqalign{
|I_1(k, t; \lambda, N)-I_2(k, t; \lambda, N)|
& \ll e^{-c_4(\log k)/\log_2k} 
{(1+e^{\lambda\kappa_t})^{2N}\log T_k
\over (e^\gamma t)^{2\kappa_t}}
\cr
& \quad
+ {E(\kappa_t)+e^{-c_4(\log k)/\log_2k}
\over N(e^\gamma t)^{2\kappa_t}}
\bigg({1+e^{\lambda\kappa_t}\over \lambda T_k}\bigg)^{2N}
\cr}
\leqno(9.6)$$
uniformly for $\lambda>0$,
$N\in \N$,
$k\ge 16$, $2\mid k$ and $t\le T(k)$,
where $T(k)$ is given by (1.10).
The implied constant depends on $\eta$ only.

\noindent{\sl Proof}.
By the definitions of $I_1$ and $I_2$, we can write
$$\eqalign{
& I_1(k, t; \lambda, N)-I_2(k, t; \lambda, N)
\cr
& \quad
= {1\over 2\pi i}\int_{\kappa_t-i\infty}^{\kappa_t+i\infty} 
\bigg(
\sum_{f\in {\rm H}_k^+(1;\eta)}\omega_f
L(1,f)^s-E(s)
\bigg)
\bigg({e^{\lambda s}-1\over \lambda s}\bigg)^{2N}
{\d s\over s(e^\gamma t)^{2s}}.
\cr}$$
In order to estimate the last integral,
we split it into two parts according to 
$|\tau|\le T_k$ or $|\tau|>T_k$.

In view of (1.18), 
it is easy to see that $\kappa_t\le T_k$ for $t\le T(k)$.
Thus we may apply (9.2) of Lemma 9.1 for 
$s=\kappa_t+i\tau$ with $|\tau|\le T_k$.
Note that $|(e^{\lambda s}-1)/(\lambda s)|\le 1+e^{\lambda\kappa_t}$
for $s=\kappa_t+i\tau$, which is easily seen 
by looking at the cases $|\lambda s|\le 1$ and $|\lambda s|>1$.
The contribution of $|\tau|\le T_k$ 
to $|I_1(k, t; \lambda, N)-I_2(k, t; \lambda, N)|$ is
$$\ll e^{-c_4(\log k)/\log_2k}
{(1+e^{\lambda\kappa_t})^{2N}\log T_k
\over (e^\gamma t)^{2\kappa_t}}.
\leqno(9.7)$$

Since $\kappa_t\le T_k$ for $t\le T(k)$, 
we can apply (9.2) of Lemma 9.1 to write, 
for $s=\kappa_t+i\tau$ with $\tau\in \R$, 
$$\eqalign{\Big|
\sum_{f\in {\rm H}_k^+(1;\eta)}\omega_f
L(1,f)^s-E(s)
\Big|
& \le \sum_{f\in {\rm H}_k^+(1;\eta)}\omega_f
L(1,f)^{\kappa_t}+E(\kappa_t)
\cr\noalign{\vskip 1mm}
& \le 2E(\kappa_t)
+ O\big(e^{-c_4(\log k)/\log_2k}\big).
\cr}$$
Thus the contribution of $|\tau|>T_k$ 
to $|I_1(k, t; \lambda, N)-I_2(k, t; \lambda, N)|$ is
$$\eqalign{
& \ll {E(\kappa_t)+e^{-c_4(\log k)/\log_2k}
\over (e^\gamma t)^{2\kappa_t}}
\int_{|\tau|\ge T_k} 
\bigg({1+e^{\lambda\kappa_t}\over \lambda|\tau|}\bigg)^{2N}
{\d\tau\over |\tau|}
\cr
& \ll {E(\kappa_t)+e^{-c_4(\log k)/\log_2k}
\over N(e^\gamma t)^{2\kappa_t}}
\bigg({1+e^{\lambda\kappa_t}\over \lambda T_k}\bigg)^{2N}.
\cr}
\leqno(9.8)$$
Combining (9.7) and (9.8) yields to the required estimate.
\hfill
$\square$

\medskip

\noindent{\it End of the proof of Theorem 2}

\smallskip

For simplicity of notation, we write
$$I_j
:=I_j(k, t; \lambda, N)
\qquad{\rm and}\qquad
I_j^+
:=I_j(k, te^{\lambda N}; \lambda, N)
\qquad(j=1, 2).$$
By using Lemma 9.2, we have
$$\eqalign{\widetilde{F}_k(t)
& \le I_1 + O\big(k^{-5/6}\big)
\cr\noalign{\vskip 0,5mm}
& = I_2 + O\big(|I_1-I_2|+k^{-5/6}\big)
\cr\noalign{\vskip 0,5mm}
& \le \Phi(te^{-\lambda N}) 
+ O\big(|I_1-I_2|+k^{-5/6}\big)
\cr\noalign{\vskip 0,5mm}
& \le \Phi(t) 
+ |\Phi(te^{-\lambda N})-\Phi(t)|
+ O\big(|I_1-I_2|+k^{-5/6}\big)
\cr}
\leqno(9.9)$$
and
$$\eqalign{\widetilde{F}_k(t)
& \ge I_1^+ + O\big(k^{-5/6}\big)
\cr\noalign{\vskip 0,5mm}
& = I_2^+ + O\big(|I_1^+-I_2^+|+k^{-5/6}\big)
\cr\noalign{\vskip 0,5mm}
& \ge \Phi(te^{\lambda N}) 
+ O\big(|I_1^+-I_2^+|+k^{-5/6}\big)
\cr\noalign{\vskip 0,5mm}
& \ge \Phi(t)-|\Phi(t)-\Phi(te^{\lambda N})|
+ O\big(|I_1^+-I_2^+|+k^{-5/6}\big).
\cr}
\leqno(9.10)$$

In view of (6.10) and Theorem 3, we have
$$\big|\Phi(t)-\Phi(te^{-\lambda N})\big|
\ll \Phi(t)
\big\{\lambda N\kappa_t(\log\kappa_t)^{1/2}
+e^{-(c_3/2)\kappa_t^\delta}\big\}$$
for $\lambda N\le e^{-t}$.
Take 
$$\lambda=e^{5A}/T_k
\qquad{\rm and}\qquad 
N=[\log_2k].
\leqno(9.11)$$
Since $T_k=e^{T(k)+{3\over 2}\log_3k+2C+\log c_5}$,
it is easy to see that
$$\lambda N\le e^{-T(k)-2C}T(k)^{-1/2}
\qquad{\rm and}\qquad
\kappa_t\asymp e^t.$$
Inserting these estimates into the preceeding inequality,
a simple calculation shows that
$$\big|\Phi(t)-\Phi(te^{-\lambda N})\big|
\le \Phi(t)
\big\{e^{t-T(k)-C}(t/T(k))^{1/2}
+O\big(e^{-c_6e^{\delta t}}\big)\big\},
\leqno(9.12)$$
provided the constant $C$ is suitably large,
where $c_6=c_6(\eta, \delta)$ is a positive constant.

Similarly by using (6.10) with $te^{\lambda N}$ in place of $t$, 
we have
$$\eqalign{
\big|\Phi(t)-\Phi(te^{\lambda N})\big|
& \ll \Phi(te^{\lambda N})
\big\{\lambda N\kappa_{te^{\lambda N}}
(\log\kappa_{te^{\lambda N}})^{1/2}
+ e^{-(c_3/2)\kappa_{te^{\lambda N}}^\delta}\big\}.
\cr}$$
Since for $t\le T(k)$ we have
$$te^{\lambda N}
= t + O\big((\log_2k)^3(\log_3k)/\log k\big)
\qquad{\rm and}\qquad
\kappa_{te^{\lambda N}}
\asymp e^{te^{\lambda N}}
\asymp e^t,$$
the preceeding estimate can be writen as
$$\eqalign{
\big|\Phi(t)-\Phi(te^{\lambda N})\big|
& \le \textstyle{1\over 4} \Phi(te^{\lambda N})
\big\{e^{t-T(k)-C}(t/T(k))^{1/2} 
+ O\big(e^{-c_6e^{\delta t}}\big)\big\}
\cr\noalign{\vskip 1mm}
& \le \textstyle{1\over 4} \Phi(t)
\big\{e^{t-T(k)-C}(t/T(k))^{1/2} 
+ O\big(e^{-c_6e^{\delta t}}\big)\big\}
\cr\noalign{\vskip 1mm}
& + \textstyle{1\over 4} \big|\Phi(t)-\Phi(te^{\lambda N})\big|
\big\{e^{t-T(k)-C}(t/T(k))^{1/2} 
+ O\big(e^{-c_6e^{\delta t}}\big)\big\},
\cr}$$
from which we deduce that
$$\big|\Phi(t)-\Phi(te^{\lambda N})\big|
\le \Phi(t)
\big\{e^{t-T(k)-C}(t/T(k))^{1/2} 
+ O\big(e^{-c_6e^{\delta t}}\big)\big\}.
\leqno(9.13)$$

By using Lemma 9.3 with $te^{\lambda N}$ in place of $t$, we have
$$\eqalign{|I_1^+-I_2^+|
& \ll e^{-c_4(\log k)/\log_2k} 
{(1+e^{\lambda\kappa_{te^{\lambda N}}})^{2N}\log T_k
\over (e^\gamma te^{\lambda N})^{2\kappa_{te^{\lambda N}}}}
\cr
& \quad
+ {E(\kappa_{te^{\lambda N}})+e^{-c_4(\log k)/\log_2k}
\over N(e^\gamma te^{\lambda N})^{2\kappa_{te^{\lambda N}}}}
\bigg({1+e^{\lambda\kappa_{te^{\lambda N}}}
\over \lambda T_k}\bigg)^{2N}.
\cr}$$
On the other hand, by using Theorem 3 and (1.25),
it is easy to see that there is a positive constant $c$ such that
$$\Phi(te^{\lambda N})
\asymp \Phi(t)
\sim {E(\kappa_t)\over 
\kappa_t\sqrt{2\pi\sigma_2}(e^\gamma t)^{2\kappa_t}} 
\gg e^{-c_8e^t/t}
\gg e^{-c_9(\log k)/[(\log_2k)^{7/2}\log_3k]}$$
for $t\le T(k)$.
Thanks to Lemma 4.5, the previous estimate can be written as 
$$|I_1^+-I_2^+|
\ll \Phi(t) 
{1\over N}
\bigg({\kappa_{te^{\lambda N}}\over 
\log\kappa_{te^{\lambda N}}}\bigg)^{1/2}
\bigg({1+e^{\lambda\kappa_{te^{\lambda N}}}
\over \lambda T_k}\bigg)^{2N}
\ll {\Phi(t)\over (\log k)^A}.
\leqno(9.14)$$

Similarly we can prove (even more easily)
$$|I_1-I_2|
\ll \Phi(t)/(\log k)^A.
\leqno(9.15)$$

Inserting (9.12) and (9.16) into (9.9)
and (9.13) and (9.15) into (9.10), we obtain
$$\widetilde{F}_k(t)
\le \Phi(t)
\big\{1+e^{t-T(k)-C}(t/T(k))^{1/2}
+ O\big(e^{-c_6e^{\delta t}}+(\log k)^{-A}\big)\big\}$$
and
$$\widetilde{F}_k(t)
\ge \Phi(t)
\big\{1-e^{t-T(k)-C}(t/T(k))^{1/2}
+ O\big(e^{-c_6e^{\delta t}}+(\log k)^{-A}\big)\big\}.
$$
This implies the first asymptotic formula of (1.13)
by taking $\eta={1\over 200}$ and $\delta={1\over 5}$.

The second can be established similarly.
This completes the proof of Theorem 2.
\hfill
$\square$

\vskip 10mm

\centerline{\pptitre \S\ 10. Proof of Theorem 1}

\bigskip

The formula (1.9) is an immediate consequence of 
Theorem 2 and (1.25).

Taking $t=T(k)$ in (1.9), we find that
$$e^{-c'_1(\log k)/\{(\log_2k)^{7/2}\log_3k\}}
\ll \widetilde{F}_k(T(k))
\ll e^{-c'_2(\log k)/\{(\log_2k)^{7/2}\log_3k\}},
\leqno(10.1)$$
where $c'_1$ and $c'_2$ are two positive constants.
Clearly (10.1) and (1.8) imply (1.11).

The related results on  $\widetilde{G}_k(t)$ and $G_k(T(k))$
can be proved similarly.
This completes the proof of Theorem 1. 
\hfill
$\square$

\vskip 8mm

\centerline{\pptitre References}

\bigskip

\item{[1]}
{\author J. Cogdell \& P. Michel},
On the complex moments of symmetric power $L$-functions at $s=1$,
{\it IMRN}. {\bf 31} (2004), 1561--1618.

\item{[2]} 
{\author D. Goldfeld, J. Hoffstein \& D. Lieman},
An effective zero-free region,
{\it Ann. of Math.} {\bf 140} (1994), 177--181, 
Appendice of [7].

\item{[3]}
{\author S.W. Graham \& C.J.  Ringrose}, 
Lower bounds for least quadratic nonresidues. 
Analytic number theory, 269--309, 
Progr. Math., 85, Birkh\"auser Boston, 1990. 

\item{[4]}
{\author A. Granville \& K. Soundararajan}, 
The distribution of values of $L(1,\chi\sb d)$,
{\it Geom. Funct. Anal.} {\bf 13} (2003), no. 5, 992--1028.

\item{[5]}
{\author L. Habsieger \& E. Royer},
$L$-functions of automorphic forms and combinatorics : Dyck paths,
{\it Ann. Inst. Fourier (Grenoble)} {\bf 54} (2004), no. 7, 
2105--2141.

\item{[6]}
{\author A. Hildebrand \& G. Tenenbaum},
On integers free of large prime factors,
{\it Trans. Amer. Math. Soc.} {\bf 296} (1986), no. 1, 265--290. 

\item{[7]}
{\author J. Hoffstein \& P. Lockhart},
Coefficients of {M}aass forms and the Siegel zero, 
{\it Ann. of Math. (2)} {\bf 140} (1994), no. 1, 161--181.

\item{[8]} 
{\author H. Iwaniec},
{\it Topics in Classical Automorphic Forms}, 
vol. 17 de {\it Graduate Studies in Mathematics}.
Amer. Math. Soc., Providence, Rhode Island, 1997.

\item{[9]}
{\author H. Iwaniec, W. Luo \& P. Sarnak},
Low lying zeros of families of $L$-functions,
{\it Inst. Hautes \'Etudes Sci. Publ. Math.} {\bf 91} (2000), 55--131.

\item{[10]}
{\author Y.-K. Lau \& J. Wu},
A density theorem on automorphic $L$-functions 
and some applications,
{\it Trans. Amer. Math. Soc.} {\bf 358} (2006), 441--472. 

\item{[11]}
{\author Y.-K. Lau \& J. Wu},
Extreme values of symmetric power $L$-functions at 1,
{\it Acta Arith.}  {\bf 126} (2007), No. 1, 57--76.

\item{[12]}
{\author W. Luo},
Values of symmetric square ${L}$-functions at $1$,
{\it J. reine angew. Math.} {\bf 506} (1999), 215--235.

\item{[13]}
{\author H.L. Montgomery \& R.C. Vaughan},
Extreme values of Dirichlet $L$-functions at $1$,
in: {\it Number theory in progress}, Vol. 2,
Zakopane-Ko\'scielisko, 1997  
(K. Gy\"ory, H. Iwaniec \& J. Urbanowicz, Eds),    
1039--1052, de Gruyter, Berlin, 1999. 

\item{[14]}
{\author E. Royer},
Statistique de la variable al\'eatoire $L(1, {\rm sym}^2f)$,
{\it Math. Ann.} {\bf 321} (2001), 667--687.

\item{[15]}
{\author E. Royer},
Interpr\'etation combinatoire des moments n\'egatifs des valeurs de
fonctions $L$ au bord de la bande critique,
{\it Ann. Sci. \'Ecole Norm. Sup. (4)} {\bf 36} (2003), 601--620

\item{[16]}
{\author E. Royer \& J. Wu},
Taille des valeurs de fonctions $L$ de carr\'es sym\'etriques 
au bord de la bande critique,
{\it Rev. Mat. Iberoamericana} {\bf 21} (2005), 263--312.

\item{[17]}
{\author E. Royer \& J. Wu},
Special values of symmetric power $L$-functions
and Hecke eigenvalues,
to appear in {\it J. Th\'eorie des Nombres de Bordeaux}.

\item{[18]} 
{\author J.-P. Serre},
{\it Abelian $\ell$-adic representation and elliptic Curves},
New York, Benjamin (1968).
Reprinted by A.K. Peters: Wellesley (1998).

\item{[19]}
{\author G. Tenenbaum}, 
{\it Introduction to analytic and probabilistic number theory},
Translated from the second French edition (1995) by C. B. Thomas,
Cambridge Studies in Advanced Mathematics {\bf 46}, 
Cambridge University Press, Cambridge, 1995. xvi+448 pp. 

\item{[20]}
{\author J. Wu},
Note on a paper by A. Granville and K. Soundararajan,
{\it J. Number Theory} {\bf 123} (2007), 329-351.

\vskip 5mm

\hskip -1mm{\author School of Mathematics and System Sciences, 
Shandong University, Jinan, Shandong 250100, China}

\hskip -1mm{\tt E-mail: jyliu@sdu.edu.cn}

\medskip

\hskip -1mm{\author Laboratoire de math\'ematiques,
UMR 6620 UBP CNRS,
Universit\'e Blaise Pascal,
F-63177 Aubi\`ere cedex,
France}

\hskip -1mm{\tt E-mail: emmanuel.royer@polytechnique.org}

\medskip

\hskip -1mm{\author Institut Elie Cartan,
UMR 7502 UHP CNRS INRIA,
Universit\'e Henri Poincar\'e (Nancy 1),
F-54506 Vand\oe uvre-l\`es-Nancy, France}

\hskip -1mm{\tt E-mail: wujie@iecn.u-nancy.fr}

\end
\bye